\documentclass[12pt]{amsart}
\usepackage[english]{babel}
 \usepackage[latin1]{inputenc}
  \usepackage[T1]{fontenc}
  \usepackage{amsmath,amssymb}
  \usepackage{amsmath}
  \usepackage{amsthm}
\usepackage{mathrsfs}
\usepackage{enumerate}
\usepackage{tikz}
\usepackage{tikz-cd}
\usetikzlibrary{shapes.geometric}
\usepackage{physics}
\usepackage[notcite, final, notref]{showkeys}
 \numberwithin{equation}{section}

\newtheorem{theorem}{Theorem}[section]

\newtheorem{lemma}[theorem]{Lemma}
\newtheorem{proposition}[theorem]{Proposition}

\newtheorem{corollary}[theorem]{Corollary}
\newtheorem{definition}[theorem]{Definition}

\theoremstyle{definition}
\newtheorem{remark}[theorem]{Remark}

\usepackage{tikz}
\usetikzlibrary{arrows,shapes,positioning,shadows,trees,matrix}
\usepackage{xcolor}

\title[]{New aspects of Bargmann transform using Touchard polynomials and hypergeometric functions}

\author[D. Alpay]{Daniel Alpay}
\address{(DA) Schmid College of Science and Technology \\
Chapman University\\
One University Drive
Orange, California 92866\\
USA}
\email{alpay@chapman.edu}

\author[A. De Martino]{Antonino De Martino}
\address{(ADM) Politecnico di
Milano\\Dipartimento di Matematica\\Via E. Bonardi, 9\\20133 Milano\\Italy}
\email{antonino.demartino@polimi.it}

\author[K. Diki]{Kamal Diki}
\address{(DK) Schmid College of Science and Technology \\
Chapman University\\
One University Drive
Orange, California 92866\\
USA}
\email{diki@chapman.edu}

\begin{document}
\maketitle
\begin{abstract}
In this paper we study the ranges of the Schwartz space $\mathcal{S}$ and its dual $\mathcal{S}'$ (space of tempered distributions) under the Segal-Bargmann transform. The characterization of these two ranges lead to interesting reproducing kernel Hilbert spaces whose reproducing kernels can be expressed respectively in terms of the Touchard polynomials and the hypergeometric functions. We investigate the main properties of some associated operators and introduce two generalized Bargmann transforms in this framework. This can be considered as a continuation of an interesting research path that Neretin started earlier in his book on Gaussian integral operators.
	\end{abstract}\mbox{}\\



\tableofcontents

\section{Introduction and preliminary results} 
The Bargmann transform introduced in \cite{Barg} is a fundamental mathematical model of quantum mechanics that maps the classical Hilbert space of square integrable functions on $\mathbb{R}^n$ onto the so-called Segal-Barmann space. This space consists of entire functions on $\mathbb{C}^n$ that are square integrable with respect to the Gaussian measure. Some authors refer to this as the \textit{boson} Fock space with $n$-degrees of freedom, see \cite{Ner}. The book \cite{Zhu} presents a comprehensive introduction to the theory of Fock spaces and related topics in mathematical analysis. For applications in the study of coherent states and other topics in quantum mechanics we refer to \cite{TAG200, A2015, Hall2013}. It is worth mentioning that Bargmann transform was extended in \cite{Asampling} to the case of polyanalytic functions using time frequency analysis tools such as the Short-time Fourier transform (STFT). A different approach using  the true polyanalytic Bargmann and polyanalytic Fock spaces was introduced by Vasilevski in \cite{Vas}. See also \cite{ACDS20222} for the case of polyanalytic functions of infinite order and \cite{FLKC2014} for the case of meta-analytic extensions. Recently new connections with the theory of superoscillations are investigated in \cite{ACDSS, STFT2023}. For recent developments in quaternionic and Clifford analysis one may consult \cite{BESC2019, CSS2016, CD2017, KA12, KA, DG2017, DKS2019, PSS2016}.\\ \\
Before presenting our main results we start by recalling some basic notions and properties of Fock spaces and Bargmann transform. For more details on this topic we refer the reader to \cite{Zhu}.  
\begin{definition}
An entire function $f: \mathbb{C} \to \mathbb{C}$ belongs to the Fock space $ \mathcal{F}(\mathbb{C})$ if
$$ \| f \|_{\mathcal{F}(\mathbb{C})}^2= \frac{1}{\pi} \int_{\mathbb{C}} | f(z)|^2 e^{-|z|^2} d \lambda(z) < \infty,$$
where $ d \lambda(z)=dx dy$ is the classical Lebesgue measure for $z=x+iy$. 
\end{definition}
The Fock space is a reproducing kernel Hilbert space, with reproducing kernel given by
$$ K(z,w)=e^{z \bar{w}}, \qquad \forall z,w \in \mathbb{C}.$$
The Fock space has the following sequential characterization
\begin{equation}
	\label{seq }
	\mathcal{F}(\mathbb{C}):= \left \{ f(z)=\sum_{n=0}^\infty z^n a_n , \quad (a_n)_{n \in \mathbb{N}_0} \subset \mathbb{C}, \quad \sum_{n=0}^\infty n!|a_n|^2  < \infty \right\}.
\end{equation}

Another fundamental tool that we use in this paper is the Segal-Bargmann transform introduced in \cite{Barg}. We use a modified version of this transform that can be found in the book \cite{Zhu}. Given a function $ \varphi \in L^{2}(\mathbb{R})$ the Segal-Bargmann transform is defined for any $z \in \mathbb{C}$ as
\begin{equation}
	\label{kernel}
	\mathcal{B}(\varphi)(z)= \int_{\mathbb{R}} A(z,x) \varphi(x) dx, \qquad A(z,x)= \left(\frac{1}{2 \pi}\right)^{\frac{1}{4}}e^{- \frac{x^2}{4}-\frac{z^2}{2}+zx}.
\end{equation}
The kernel of the Segal Bargmann transform can be written in terms of generating function of the normalized Hermite functions. Precisely,  the Hermite functions are defined by
\begin{equation}
\label{harmite}
h_n(x)=e^{- \frac{x^2}{4}}H_n(x)=(-1)^n e^{\frac{x^2}{2}}\frac{d^n}{d x^n} \left( e^{-\frac{x^2}{2}}\right),
\end{equation}
where $H_n(x)$ are the probabilistic Hermite polynomials. Hence the normalized Hermite functions are defined as

$$ \xi_n(x)= \frac{h_n(x)}{\| h_n\|_{L^2}}= \frac{(-1)^n}{(2\pi)^{\frac{1}{4}} \sqrt{n!}}e^{\frac{x^2}{2}} \frac{d^n}{d x^n} \left(e^{-\frac{x^2}{2}}\right).$$
Therefore we can write the kernel of the Bargmann transform as
\begin{equation}
\label{genkern}
A(z,x)=\frac{1}{(2 \pi)^{\frac{1}{4}}} \sum_{n=0}^{\infty} \frac{z^n h_n(x)}{n!}.
\end{equation}
It is well-known that $\mathcal{B}$ is an unitary transform mapping $L^2(\mathbb{R})$ onto the Fock space $\mathcal{F}(\mathbb{C})$. The inverse of the Bargmann transform is denoted by $ \mathcal{B}^{-1}$. Moreover, the Bargmann transform applied to the normalized Hermite functions gives back complex monomials:
\begin{equation}
\label{prop1}
\mathcal{B}(\xi_n)(z)= \frac{z^n}{\sqrt{n!}}, \qquad n \in \mathbb{N}_0.
\end{equation}

 The main purpose of this paper is to continue exploring this fascinating research topic and investigate new aspects of this transform involving relations with the Schwartz space and its dual space of tempered distributions. Indeed, inspired by the discussion in \cite[Proposition 2.5 and Proposition 2.6]{Ner} we aim to give a complete characterization of the range of the Schwartz space $\mathcal{S}$ and its dual $\mathcal{S}'$ under the Bargmann transform. We will develop a detailed  study of these spaces. In literature the Bargmann transform is studied on the Hilbert space $L^2$ that contains the Schwartz space $\mathcal{S}$ and which is contained in the space of tempered distributions $\mathcal{S}'$. Namely, we have 
\begin{equation}
\label{star1}
\mathcal{S}\subset L^2\subset \mathcal{S}'.
\end{equation}
It is well-known that the range of the space $L^2$ under the Bargmann transform is the Fock space. It is natural to wonder if there exists a characterization of the ranges of the Schwartz space $\mathcal{S}$ and the space of tempered distributions $\mathcal{S}'$ under the Bargmann transform. In this paper we tackle this problem and we find out that the action of the Bargmann transform on the spaces $\mathcal{S}$ and $\mathcal{S}'$ lead to interesting reproducing kernel Hilbert spaces. Precisely, in this paper we introduce and study the following spaces
$$	\mathcal{H}_p(\mathbb{C}):= \left\{ f(z)=\sum_{n=0}^\infty z^n a_n ,\,  \quad (a_n)_{n \in \mathbb{N}_0} \subseteq \mathbb{C}, \quad \sum_{n=0}^\infty \frac{|a_n|^2 n!}{(n+1)^{2p}}   <\infty \right\},
$$
$$	\mathcal{F}_p(\mathbb{C}):= \left\{ f(z)=\sum_{n=0}^\infty z^n a_n ,\,  \quad (a_n)_{n \in \mathbb{N}_0} \subseteq \mathbb{C}, \quad \sum_{n=0}^\infty |a_n|^2 n!(n+1)^{2p} <\infty \right\},$$
for $p>0$.
The reproducing kernels of the spaces $\mathcal{H}_p(\mathbb{C})$ and $ \mathcal{F}_p(\mathbb{C})$ can be expressed in terms of the Touchard polynomials and the hypergeometric functions. The Touchard polynomials were introduced first in \cite{T} as polynomials having Stirling numbers of the second kind as coefficients. They share same interesting properties with other well-known families of orthogonal  polynomials such as: Rodrigues formula, generating functions, orthogonality conditions, etc. \\ \\ Associated to the spaces $ \mathcal{H}_p(\mathbb{C})$ and $ \mathcal{F}_p(\mathbb{C})$ we introduce and study two generalized Bargmann transforms that are constructed by taking the generating functions associated to some orthonormal basis against the Hermite functions. It is also possible to think about the Bargmann construction here by acting on the classical Segal-Bargmann with some specific operators combining powers of the momentum and annihilation operators. The technique that we use to achieve these objectives consists of computing higher derivatives of the classical Segal-Bargmann kernel and taking their linear combinations with suitable coefficients involving the Stirling numbers of the second kind. In particular, while developing the results of this paper we computed explicitly in Proposition \ref{der} the higher derivatives of the classical Segal-Bargmann kernel $A(z,x)$, which is given by

$$
\frac{d^k}{dz^k}A(z,x)=H_k(x-z)A(z,x),
$$

where $k\geq 0$ and $H_k(z)$ denotes the holomorphic extension of the Hermite polynomials, see \cite{S, V}. 
\\ \\ 

\emph{Plan of the paper:} in Section 2 we construct the new reproducing kernel Hilbert space $\mathcal{H}_p$ that are built upon the characterization of the range of the space of tempered distributions under the Bargmann transform. In particular, we obtain a new reproducing kernel function whose expression involves the Touchard polynomials. We also study some operators associated to the $\mathcal{H}_p$ space and investigate a generalized Bargmann transform based on this space. In Section 3 we study the reproducing kernel Hilbert space $\mathcal{F}_p$ that characterizes the range of the Schwartz space under the Bargmann transform. Finally, in Section 4 we present a connection between the reproducing kernel Hilbert space $\mathcal{H}_p$ and $\mathcal{F}_p$, through the Fock space.

\section{The space $\mathcal{H}_p$ and associated operators}

\subsection{Construction and basic properties of the $\mathcal{H}_p$-space}
It is well-known that the the Schwartz space can be seen as intersection of weighted $\ell^2$ spaces, whereas the union of weighted $\ell^2$ space give the space of tempered distribution, see \cite{Sim}. Precisely we have

$$ \bigcap_{p>0} \ell^2(\mathbb{N}, (n+1)^{2p})= \mathcal{S}(\mathbb{R}) \subset L^2(\mathbb{R}) \subset \mathcal{S}'(\mathbb{R})= \bigcup_{p>0} \ell^{2}(\mathbb{N}, (n+1)^{-2p}).$$
\begin{lemma}
\label{star}
For a fixed $p > 0$ we have
$$ \mathcal{B} \left(\ell^2(\mathbb{N}, (n+1)^{-2p})\right)=\mathcal{H}_p(\mathbb{C}),$$
where
\begin{equation}
	\label{space}
	\mathcal{H}_p(\mathbb{C}):= \left\{ f(z)=\sum_{n=0}^\infty z^n a_n ,\,  \quad (a_n)_{n \in \mathbb{N}_0} \subseteq \mathbb{C}, \quad \sum_{n=0}^\infty \frac{|a_n|^2 n!}{(n+1)^{2p}}   <\infty \right\}.
\end{equation}
\end{lemma}
\begin{proof}
To show the result we have to prove that
$$ \mathcal{B}(S_p)= \mathcal{H}_p(\mathbb{C}), $$
where
$$S_p:=\left\{\sum_{n=0}^\infty \xi_n a_n, \quad (a_n)_{n \in \mathbb{N}_0} \subseteq \mathbb{C}, \quad \sum_{n=1}^\infty (n+1)^{-2p} |a_n|^2< \infty\right\}.$$
A function $g$ belongs to $\mathcal{B}(S_p)$ if and only if there exists $f \in S_p$ such that $g(z)= \mathcal{B}(f)(z)$. From the action of the Bargmann transform applied to the nominalized Hermite functions, see \eqref{prop1}, we get 
$$ \mathcal{B}(f)(z)=\sum_{n =1}^\infty \frac{z^n}{\sqrt{n!}} a_n.$$
We set $b_n:= \frac{a_n}{\sqrt{n!}}$. Then we have
$$ \sum_{n=1}^\infty |b_n|^2  \frac{n!}{(n+1)^{2p}}= \sum_{n=1}^\infty |a_n|^2 (n+1)^{-2p}< \infty.$$
Therefore we get $ \mathcal{B}(f)(z) \in \mathcal{H}_p(\mathbb{C})$
\end{proof}
Thanks to \eqref{prop1} the Bargmann transform can be defined on the whole of $  \mathcal S^\prime(\mathbb{R})$. This extension, still denoted by the letter
$\mathcal B$, is defined by
\[
\mathcal B((a_n)_{n\in\mathbb N_0})=\sum_{n=0}^\infty\frac{z^n}{\sqrt{n!}}a_n,
\]
or, formally, by
\[
\mathcal B\left(\sum_{n=0}^\infty \xi_n a_n\right)=\sum_{n=0}^\infty\frac{z^n}{\sqrt{n!}}a_n.
\]
We have the following result, characterizing the image of $\mathcal S^\prime$ under $\mathcal B$.

\begin{theorem}
	$\mathcal B$ is bicontinuous from $\mathcal{S'}$ onto $\mathcal B(\mathcal S^\prime)=\bigcup_{p \geq 1} \mathcal{H}_p(\mathbb{C})$.
\end{theorem}
\begin{proof}
We denote by $i_p$ the injection map from   $\ell^2(\mathbb N, (n+1)^{-2p})$ into
	\[
	\mathcal S^\prime=\cup_{q=0}^\infty \ell^2(\mathbb N, (n+1)^{-2q}).
	\]
	The inductive topology on the latter is by definition the finest topology for which all the injections $i_p$ are continuous,
	see e.g. \cite[Definition 1 p. 136]{groth}.
	A similar remark holds for the image of $\mathcal S^\prime$ under $\mathcal B$.\smallskip
	
	By the characterization of continuous maps with domain a space, with an inductive topology and domain a locally convex space
	(see \cite[\S 6.1 p. 54]{schaefer}), we have to show that for every $p$ the map
	$\mathcal B\circ i_p$ is continuous from $\ell^2(\mathbb N, (n+1)^{-2p})$ into  $\mathcal B(\mathcal S^\prime)$.
	But the restriction $\mathbf{B}_p$ of the Bargmann map $\mathcal B$ to $\ell^2(\mathbb N, (n+1)^{-2p})$ has range
	$\mathcal B(\ell^2(\mathbb N, (n+1)^{-2p})=\mathcal H_p(\mathbb C)$, see Lemma \ref{star}, and is continuous since it is isometric by construction.
	Furthermore, the inclusion map $I_p$ from $\mathcal B(\ell^2(\mathbb N, (n+1)^{-2p}))$
	into $\mathcal B(\mathcal S^\prime)$ is continuous by the definition of the inductive topology on $\mathcal B(\mathcal S^\prime)$.
	But for $a\in(\ell^2(\mathbb N, (n+1)^{-2p}))$ we have
	\[
	i_p(a)=a\quad{\rm and}\quad \mathbf{B}_p(a)= \mathcal B(a).
	\]
	Hence,
	\[
	\mathcal B\circ i_p(a)=\mathcal B(a)\quad{\rm and}\quad   I_p\circ \mathbf{B}_p(a)=\mathbf{B}_p(a)= \mathcal B(a).
	\]
	Thus
	\[
	\mathcal B\circ i_p=I_p\circ \mathbf{B}_p
	\]
	so that $\mathcal B\circ i_p$ is continuous since $I_p\circ \mathbf{B}_p$ is continuous. Thus $\mathcal B$ is continuous by definition of the inductive topology.\smallskip
	
	The proof for $\mathcal B^{-1}$ is similar.  
\end{proof}
Now, we aim to discuss the main properties of the space $ \mathcal{H}_p(\mathbb{C})$. This space takes an importance  because it constitutes the building block of the range  of the space $ \mathcal{S}'$ under the Bargmann transform. First we give a rigorous definition of this space.
\begin{definition}
Let $p \in \mathbb{N}$. We define the Hilbert space $ \mathcal{H}_p(\mathbb{C})$ as the set of all entire functions endowed with the inner product
\begin{equation}
\label{inner0}
\langle f, g \rangle_{\mathcal{H}_p(\mathbb{C})}= \sum_{n=0}^\infty \frac{n!}{(n+1)^{2p}} a_n \overline{b_n},
\end{equation}
where $ f(z)= \sum_{n=0}^\infty z^n a_n$ and $g(z)= \sum_{n=0}^\infty b_n z^n$, with $ (a_n)_{n \in \mathbb{N}_0}$, $(b_n)_{n \in \mathbb{N}_0} \subseteq \mathbb{C}$.
\end{definition}
The inner product defined in \eqref{inner0} induces the following norm
$$ \| f \|_{\mathcal{H}_p(\mathbb{C})}^2= \sum_{n=0}^\infty \frac{n!}{(n+1)^{2p}} | a_n|^2, \qquad f \in \mathcal{H}_p(\mathbb{C}).$$

\begin{remark}
For $p=0$ the space $ \mathcal{H}_p(\mathbb{C})$ coincides with the classical Fock space, see \eqref{seq }.
\end{remark}
As a consequence of \eqref{inner0} the inner product of two different monomials is given by the following formula
\begin{equation}
\label{inner}
\langle z^n, z^m \rangle_{\mathcal{H}_p(\mathbb{C})}:= \frac{n!}{(n+1)^{2p}} \delta_{n,m}, \qquad n,m \geq 0.
\end{equation}
Now, we consider the following function
\begin{equation}
	\label{kern}
	K_p(z,w)=\sum_{n=0}^\infty \frac{z^n \bar{w}^n (n+1)^{2p}}{n!}.
\end{equation}
\begin{proposition}
\label{necc}
Let $w \in \mathbb{C}$ and $p \in \mathbb{N}$. 
\begin{itemize}
\item[i)] The function $K_p(z,w)$ belongs to $ \mathcal{H}_p(\mathbb{C})$.
\item[ii)] The evaluation map $e_w(f):=f(w)$ is a continuous linear functional on $ \mathcal{H}_p(\mathbb{C})$. Furthermore, we have
$$|e_w(f)| \leq \| f\|_{\mathcal{H}_p(\mathbb{C})}\| K_p\|_{\mathcal{H}_p(\mathbb{C})}, \qquad \forall f \in \mathcal{H}_p(\mathbb{C}).$$
\end{itemize}
\end{proposition}
\begin{proof}
\begin{itemize}
\item[i)]  We set $ \alpha(w):= \frac{\bar{w}^n (n+1)^{2p}}{n!}$ and thus we can write
$$ K_p(z,w)= \sum_{n=0}^{\infty} z^n \alpha(w).$$
Hence we have
\begin{equation}
	\label{norm}
\| K_p \|_{\mathcal{H}_p(\mathbb{C})}=\sum_{n=0}^{\infty} \frac{n!}{(n+1)^{2p}}|\alpha(w)|^2=\sum_{n=0}^{\infty}  \frac{(n+1)^{2p}}{n!} | w|^n < \infty. 
\end{equation}
\item[ii)] Let $(a_n)_{n \in \mathbb{N}_0} \subseteq \mathbb{C}$. By the Schwartz inequality and \eqref{norm} we have
\begin{eqnarray*}
|e_w(f)| & \leq & \sum_{n=0}^{\infty} |w|^n |a_n|\\
& =& \sum_{n=0}^{\infty} \frac{\sqrt{n!}}{(n+1)^{p}} \frac{(n+1)^p}{\sqrt{n!}} |w|^n |a_n|\\
& \leq&  \left( \sum_{n=0}^{\infty} \frac{(n+1)^{2p}}{n!} |a_n|^2\right)^{\frac{1}{2}} \left( \frac{n!}{(n+1)^{2p}} |w|^{2n}\right)^{\frac{1}{2}}\\
& \leq & \| f\|_{\mathcal{H}_p(\mathbb{C})}\| K_p\|_{\mathcal{H}_p(\mathbb{C})}.
\end{eqnarray*}
\end{itemize}
\end{proof}
The previous result guarantees the existence of a reproducing kernel for the space $\mathcal{H}_p$. 

\begin{theorem}
\label{kern2}
Let $p \in \mathbb{N}$, the reproducing kernel of the space $\mathcal{H}_p(\mathbb{C})$ is given by $K_p(z,w)$.

\end{theorem}
\begin{proof}
Let $f(z)= \sum_{n=0}^\infty z^n a_n$, with $(a_n)_{n \geq 0}\subseteq \mathbb{C}$, being in the space $ \mathcal{H}_p(\mathbb{C})$. We set $b_n:= \frac{ \bar{w}^n (n+1)^{2p}}{n!}$. By formula \eqref{inner0} we have
$$ \langle f, K_p \rangle_{\mathcal{H}_p(\mathbb{C})}=\sum_{n=0}^\infty \frac{n!}{(n+1)^{2p}} a_n \overline{b}_n=\sum_{n=0}^{\infty} a_n w^n=f(w).$$

This proves the reproducing kernel property.
\end{proof}
Now, we establish a relation between the reproducing kernel introduced in \eqref{kern} and a special class of polynomials: the Touchard polynomials. These were first introduced in \cite{T} and are defined as
$$ T_n(x):= \sum_{k=0}^n S(n,k) x^k,$$
where $S(n,k)$ are the Stirling numbers of the second kind and can be calculated using the following explicit formula  
\begin{equation}
\label{Stirling}
S(n,k):=\frac{1}{k!} \sum_{i=0}^k (-1)^i \binom{k}{i} (k-i)^n.
\end{equation}
Few examples of Touchard polynomials are given by
$$ T_0(x)=1,$$
$$
 T_1(x)=x,
$$
$$ T_2(x)=x^2+x,$$
$$ T_3(x)=x^3+3x^2+x,$$
$$ T_4(x)=x^4+6x^3+7x^2+x,$$
$$ T_5(x)= x^5+10x^4+25x^3+15x^2+x.$$
The generating function for the polynomials $T_n(x)$ is given by 
\begin{equation}
\label{proper}
T_n(x)=e^{-x} \sum_{k=0}^\infty \frac{k^n x^k}{k!}.
\end{equation}

\begin{remark}
	It is clear that the Touchard polynomials can be easily extended in complex variables.  
\end{remark}

\begin{proposition}
\label{rel}
Let $p \in \mathbb{N}$, then the reproducing kernel of the space $\mathcal{H}_p(\mathbb{C})$ can be written as
\begin{equation}
\label{repkernel}
K_p(z \bar{w})=\frac{e^{z \bar{w}}}{z \bar{w}} T_{2p+1}(z\bar{w}).
\end{equation}
\end{proposition}
\begin{proof}
By formula \eqref{proper} and Theorem \ref{kern2} we have that
\begin{eqnarray*}
\frac{e^{z \bar{w}}}{z \bar{w}} T_{2p+1}(z\bar{w})&=& \frac{1}{z \bar{w}} \sum_{k=0}^\infty \frac{k^{2p+1} z^k \bar{w}^k}{k!}\\
&=&   \sum_{k=1}^\infty \frac{k^{2p} z^{k-1} \bar{w}^{k-1}}{(k-1)!}\\
&=& \sum_{k=0}^\infty \frac{(k+1)^{2p} z^k \bar{w}^k}{k!}\\
&=& K_p(z \bar{w}).
\end{eqnarray*}
\end{proof}
\begin{remark}
The Touchard polynomials $T_{2p+1}(z \bar{w})$, for $p \geq 1$, can be always factorized in terms of $z \bar{w}$. Therefore formula \eqref{repkernel} does not present any singularity.
\end{remark}
\begin{remark}
If we consider $p=0$ in \eqref{repkernel}, since $T_{1}(z \bar{w})= z \bar{w}$, we get back the reproducing kernel of the classical Fock space. 
\end{remark}

Now we tackle the problem to find the existence of a geometric description of the space $ \mathcal{H}_p(\mathbb{C})$. This is related to study the one-dimensional (full)-moment problems:
\newline
\newline
\emph{Given a real sequence $s=(s_n)_{n \in \mathbb{N}_0}$ and a closed subset $K$ of $ \mathbb{R}$, the $K$-moment problems states: When does there exist a Radon measure $\mu$ on $\mathbb{R}$ supported on $K$ such that $ s_n= \int_{\mathbb{R}} x^n d \mu(x)$ for all $n \in \mathbb{N}_0$?}
\newline
\newline
For $K=[0, + \infty)$ the solution of the previous problem is given by the so-called Stieltjes theorem, see \cite[Thm. 3.12]{KS}. In order to state this result we need to fix some notations. We denote by $Es$ the shifted sequence given by
$$ (Es)_{n}=s_{n+1}, \qquad n \in \mathbb{N}_0.$$
Furthermore, we denote the Hankel matrix $\mathsf{H}_n(x)$ by
$$ \mathsf{H}_n(x)= \begin{pmatrix}
s_0 && s_1 && s_2 &&...&& s_n\\
s_1 && s_2 && s_3 && ...&&s_{n+1}\\
s_2 && s_3 && s_4 &&...&& s_{n+2}\\
...&&...&&...&&...&&\\
s_n && s_{n+1} && s_{n+2} &&...&& s_{2n}
\end{pmatrix}
$$

\begin{theorem}[Solution of the Stieltjes moment problem]
\label{hamb}
A real sequence $s=(s_n)_{n \in \mathbb{N}_0}$ is a Stieltjes sequence, that is, there exists a Radon measure $ \mu$ on $[0, \infty)$  such that
\begin{equation}
\label{star3}
s_n= \int_0^\infty x^n d \mu(x), \quad  \hbox{for} \quad n \in \mathbb{N}_0,
\end{equation}
if and only if all the Hankel matrices $\mathsf{H}_n(s)$ and $\mathsf{H}_n(Es)$, $n \in \mathbb{N}_0$, are positive semidefinite.
\end{theorem}
In \cite[Cor. 4.10]{KS} the author proved a condition that the Stieltjes moment problem has to satisfy in order to have a unique solution.

\begin{proposition}
\label{sequence}
Let $s=(s_n)_{n \in \mathbb{N}_0}$ being a Stieltjes sequence. Then if there exists a constant $M >0$ such that
$$
s_n \leq M^n (2n)! \qquad \hbox{for} \quad n \in \mathbb{N},
$$
then the sequence $s$ has only one representing measure on $[0, \infty)$.
\end{proposition}

To show the existence of a geometric description for the space $ \mathcal{H}_p(\mathbb{C})$ we need also the following generic result of changing variables for an integral, see \cite[Prop. 4.10]{RY}.
\begin{theorem}
\label{change1}
If $M$ is a continuous, non decreasing function on the interval $[a,b]$, then for $g$ being a non negative Borel function on $[M(a), M(b)]$ we have
$$\int_{M(b)}^{M(a)} g(\rho) dN(\rho)=\int_{a}^b g(M(\rho)) d N(M(\rho)).$$
\end{theorem}
Another fundamental tool that we use is the Hadamard (also called Schur) product between matrices. Given two matrices $A \in \mathbb{C}^{m \times n}$ and $B \in \mathbb{C}^{m \times n}$ the Hadamard product is denoted by $A \odot B$ and its elements are given by
$$ (A \odot B)_{ij}:= (A)_{ij} (B)_{ij}
$$
The peculiarity of this product of matrices is that the product of two  positive-semidefinite matrices is still  positive-semidefinite. This result is known as Schur product theorem.
\\Now, we have all the tools to prove the existence of a geometric description for the space $ \mathcal{H}_p(\mathbb{C})$.

\begin{theorem}
There exists a unique radial measure $d \sigma(x,y)$ such that
$$ \frac{1}{2 \pi} \int_{\mathbb{C}} z^n \bar{z}^m d \sigma(x,y)= \frac{n!}{(n+1)^{2p}} \delta_{n,m}, \qquad z=x+iy.$$
\end{theorem}

\begin{proof}
We start by passing in polar coordinates in the integral $\int_{\mathbb{C}} z^n \bar{z}^m d \sigma(x,y)$. Let $z= \rho e^{i \theta}$, where $ \rho \in [0, \infty)$ and $ \theta \in [0, 2 \pi]$, then we have
\begin{eqnarray*}
\frac{1}{2 \pi} \int_{\mathbb{C}} z^n \bar{z}^m d \sigma(x,y) &=& \frac{1}{2 \pi} \int_{0}^\infty \int_{0}^{2 \pi} \rho^{n+m} e^{i(n-m) \theta} d \theta d \mu (\rho)\\
&=& \delta_{n,m} \int_{0}^\infty \rho^{2n} \mu(\rho).
\end{eqnarray*}
This means that the problem is reduced to find a measure $ \mu(\rho)$ such that
\begin{equation}
\label{prob2}
\int_{0}^\infty \rho^{2n}  d \mu (\rho)= \frac{n!}{(n+1)^{2p}}.
\end{equation}
Now, we use Theorem \ref{change1} with $M(\rho)=\rho^2$, $ g(\rho)=\rho^n$, $dN= d \nu$, $a=0$ and $b=k \in \mathbb{R}$, we get
$$\int_{0}^k \rho^{2n} d \nu(\rho^2) =\int_{0}^{k^2} \rho^n d \nu(\rho) .$$
By sending $k$ to infinity we get
$$\int_{0}^\infty \rho^{2n} d \nu(\rho^2)= \int_{0}^{\infty} \rho^n d \nu(\rho).$$
If we set $d \mu(\rho):= d \nu (\rho^2)$ we get
$$\int_{0}^\infty \rho^{2n} d \mu(\rho)= \int_{0}^{\infty} \rho^n d \nu(\rho).$$
By the expression \eqref{prob2} we have to find a measure $ d \nu(\rho)$ such that
\begin{equation}
\label{hamb2}
\int_{0}^\infty \rho^n d \nu(\rho)= \frac{n!}{(n+1)^{2p}}.
\end{equation}
Now, the problem is reduced to the so-called Stieltjes-moment problem.
\\We set $s_n:= \frac{n!}{(n+1)^{2p}}$. In order to show that $(s_n)_{n \in \mathbb{N}_0}$ is a Stieltjes sequence we have to use Theorem \ref{hamb}. Hence we have to show that
$$ T_N:= \left( \frac{(n+m)!}{(n+m+1)^{2p}}\right)_{n,m=0}^N \quad \hbox{and} \quad S_N:= \left(\frac{(n+m-1)!}{(n+m)^{2p}}\right)_{m,n=1}^N$$
are positive semidefinite matrices.
We observe that we can write the matrix $T_N$ as the Hadamard product of two matrices 
$$ \left((n+m)!\right)_{n,m=0}^N \odot \left(\frac{1}{(m+m+1)^{2p}}\right)_{n,m=0}^N.$$
It is clear that the matrix
$$ \left(\frac{1}{(m+m+1)^{2p}}\right)_{n,m=0}^N$$
is positive semidefinite. Now, we observe that
$$ \left((n+m)!\right)_{m,n=0}^N= \hbox{Diag}(m!) \left(\frac{(n+m)!}{n! m!}\right)_{n,m=0}^N \hbox{Diag}(n!),$$
where we denote by $\hbox{Diag}(m!)$ and $\hbox{Diag}(n!)$ the matrices that have entries $m!$ and $n!$ on the diagonal, respectively.
By the Chu-Vandermonde formula 
$$ \binom{u+k}{k}= \sum_{\ell=0}^{k} \binom{k}{\ell} \binom{u}{\ell}$$
we get that
$$ \left(\frac{(n+m)!}{n! m!}\right)_{n,m=0}^N$$
is positive semidefinite. This implies that the matrix $T_N$ is positive-semidefinite. To show that also the matrix $S_N$ is positive definite we observe that we can write
$$S_N= \left((n+m-1)!\right)_{m,n=1}^N \odot \left( \frac{1}{(n+m)^{2p}}\right)_{n,m=1}^N.$$
It is clear that $\left( \frac{1}{(n+m)^{2p}}\right)_{n,m=1}^N$ is positive-semidefinite. Now, we show that $\left((n+m-1)!\right)_{m,n=1}^N$ is also positive semidefinite. This follows from the fact that
$$ (n+m-1)!= \Gamma(n+m-2)= \int_{0}^\infty t^{n+m-2} e^{-t} dt.$$
Therefore we get that also the matrix $S_N$ is positive-semidefinite. Hence by Theorem \ref{hamb} there exists a measure $ d \nu(\rho)$ that admits a representation like the one in \eqref{hamb2}. 
\\The unicity of the solution follows by Corollary \ref{sequence}. Indeed by using the Stirling's approximation it is not difficult to show that
$$ s_n \leq M^{n}(2n!), \qquad \forall n \in \mathbb{N}_0.$$
This concludes the proof. 
\end{proof}

\subsection{Some operators on $ \mathcal{H}_p$-space}
 The classical Fock space is the unique Hilbert space of power series defined in a neighbourhood of the origin for which the annihilation operator $\partial$ and the momentum operator $M_z$, namely
 $$ M_z f(z):=z f(z), \qquad \partial f(z):= \frac{d}{dz} f(z).$$
 are closed defined on the span of polynomials and adjoint to each other, i.e.
 $$\partial^{*}=M_z.$$
 
 Moreover the annihilation and the momentum operators satisfy the classical commutation rule
\begin{equation}
\label{comm}
[M_z, \partial]=  M_z\partial-\partial M_z =-\mathcal{I},
\end{equation}
 where $\mathcal{I}$ is the identity operator.
 \\In this context the backward-shift operator $R_0$ is defined as
 $$ R_0= \begin{cases}
 	\frac{f(z)-f(0)}{z}, \qquad z \neq 0\\
 	f'(0), \qquad \quad z=0,
 \end{cases}
 $$
 for analytic functions in a neighbourhood of the origin. Another operator that play a role in this context is the integration operator $I$, introduced in \cite{ACK} and defined by
 $$(If)(z)=\int_{[0,z]}f(s)ds.$$
 In \cite[Lemma 2.1]{ACK}the authors proved that the integration operator is the adjoint of the backward shift operator,i.e.
 $$R_0^{*}=I.$$
The commutator of $R_0$ and $I$ is equal on the monomials to
 \begin{equation}
 \label{comm5}
[R_0, I](z^n)= \begin{cases}
1, \qquad n=0\\
-\frac{z^n}{n(n+1)}, \quad n=1,2,...
\end{cases}
 \end{equation}
We denote by $D_0$ the (formal) diagonal operator defined by
\begin{equation}
\label{mat}
D_0:= \hbox{diag} \left(1,- \frac{1}{2},- \frac{1}{6},..., -\frac{1}{n(n+1)},...\right).
\end{equation}
Hence we can write
\begin{equation}
\label{com}
[R_0, I]=D_0.
\end{equation}

The goal of this section is to study the adjoint of the operators $R_0$, $\partial$, $M_z$ and $I$ for the space $ \mathcal{H}_p(\mathbb{C})$.

\begin{proposition}
\label{action2}
Let $p \in \mathbb{N}$. The adjoint operator of the backward-shift operator in $ \mathcal{H}_p(\mathbb{C})$ is given by
\begin{equation}
\label{action1}
R_0^* z^n= \frac{(n+2)^{2p}}{(n+1)^{2p+1}} z^{n+1}.
\end{equation}

\end{proposition}
\begin{proof}
We have to show that $ \langle R_0 z^n,  z^m \rangle_{\mathcal{H}_p(\mathbb{C})}=\langle  z^n,  R_0^* z^m \rangle_{\mathcal{H}_p(\mathbb{C})}$ for $m$, $n \in \mathbb{N}_0$.
If we suppose $n \geq 1$, by formula \eqref{inner} we have 
$$\langle  R_0 z^n, z^m \rangle_{\mathcal{H}_p(\mathbb{C})}=\langle  z^{n-1}, z^{m} \rangle_{\mathcal{H}_p(\mathbb{C})}=\frac{(n-1)!}{n^{2p}} \delta_{n-1, m},$$
and
$$
	\langle  z^n, R_0^*z^m \rangle_{\mathcal{H}_p(\mathbb{C})}=\frac{(m+2)^{2p}}{(m+1)^{2p+1}} \langle  z^{n}, z^{m+1} \rangle_{\mathcal{H}_p(\mathbb{C})}= \frac{(m+2)^{2p}}{(m+1)^{2p+1}} \frac{n!}{(n+1)^{2p}}\delta_{n,m+1}= \frac{(n-1)!}{n^{2p}} \delta_{n,m+1}.
$$
This proves the result for $n \geq 1$.	For $n=0$ the result is trivial, indeed we have
$$\langle R_0 1, z^m \rangle_{\mathcal{H}_p(\mathbb{C})}=0 \qquad \hbox{and} \qquad \langle  1,R_0^* z^m \rangle_{\mathcal{H}_p(\mathbb{C})}= \frac{(m+2)^{2p}}{(m+1)^{2p}} \langle  1, z^{m+1} \rangle_{\mathcal{H}_p(\mathbb{C})}=0.$$
\end{proof}

\begin{proposition}
\label{action3}
Let $p \in \mathbb{N}$. The adjoint operator of the annihilation operator in $ \mathcal{H}_p(\mathbb{C})$ is given by
\begin{equation}
\label{action0}
\partial^{*} z^n= \frac{(n+2)^{2p}}{(n+1)^{2p}} z^{n+1}.
\end{equation}

\end{proposition}
\begin{proof}
For $m$, $n \in \mathbb{N}_0$ we have to show that $\langle \partial z^n , z^m \rangle_{\mathcal{H}_p(\mathbb{C})}= \langle  z^n , \partial^* z^m \rangle_{\mathcal{H}_p(\mathbb{C})}$. For $n \geq 1$ and by formula \eqref{inner} we get

$$
 \langle  \partial z^n ,  z^m \rangle_{\mathcal{H}_p(\mathbb{C})}= n \langle  z^{n-1} ,z^{m} \rangle_{\mathcal{H}_p(\mathbb{C})}=n \frac{(n-1)!}{n^{2p}} \delta_{n-1, m}=\frac{n!}{n^{2p}} \delta_{n-1, m}.
$$
By formula \eqref{action0} we have
$$\langle  z^n , \partial^* z^m \rangle_{\mathcal{H}_p(\mathbb{C})}=  \frac{(m+2)^{2p}}{(m+1)^{2p}} \langle z^{n} , z^{m+1} \rangle_{\mathcal{H}_p(\mathbb{C})}= \frac{(m+2)^{2p}}{(m+1)^{2p}}\frac{n!}{(n+1)^{2p}} \delta_{n,m+1}=\frac{n!}{n^{2p}}\delta_{n,m+1}.$$
Hence we have the result. For $n=0$ the equality is trivial. Indeed we have
$$ \langle \partial 1 , z^m \rangle_{\mathcal{H}_p(\mathbb{C})}= 0 \qquad \hbox{and} \qquad \langle 1, \partial^{*} z^m \rangle_{\mathcal{H}_p(\mathbb{C})}= \frac{(m+2)^{2p}}{(m+1)^{2p}} \langle 1, z^{m+1}\rangle_{\mathcal{H}_p(\mathbb{C})}=0.$$

\end{proof}

\begin{proposition}
\label{action4}
Let $p \in \mathbb{N}$. Then adjoint momentum operator in $ \mathcal{H}_p(\mathbb{C})$ is given by
$$ M_z^* z^n= \frac{n^{2p+1}}{(n+1)^{2p}} z^{n-1}, \qquad M_z^*(1)=0$$
\end{proposition}
\begin{proof}
For $m$, $n \in \mathbb{N}_0$ we have to show $\langle M_z z^n, z^m\rangle_{\mathcal{H}_p(\mathbb{C})}=\langle  z^n, M_z^* z^m\rangle_{\mathcal{H}_p(\mathbb{C})}$. For $m \geq 1$ and by formula \eqref{inner} we have
$$ \langle  M_z z^n,  z^m\rangle_{\mathcal{H}_p(\mathbb{C})}=\langle  z^{n+1},  z^{m}\rangle_{\mathcal{H}_p(\mathbb{C})}= \frac{(n+1)!}{(n+2)^{2p}} \delta_{n+1, m},$$
while
$$  \langle  z^n,M_z^* z^m\rangle_{\mathcal{H}_p(\mathbb{C})}=  \frac{m^{2p+1}}{(m+1)^{2p}} \langle z^{n}, z^{m-1} \rangle_{\mathcal{H}_p(\mathbb{C})}=\frac{m^{2p+1}}{(m+1)^{2p}} \frac{n!}{(n+1)^{2p}} \delta_{n,m-1}=\frac{(n+1)!}{(n+2)^{2p}} \delta_{n,m-1}.$$
This proves the result for $m \geq 1$. For $m=0$ the result is trivial, indeed we have
$$\langle  M_z z^n, 1\rangle_{\mathcal{H}_p(\mathbb{C})}=\langle  z^{n+1},1\rangle_{\mathcal{H}_p(\mathbb{C})}=0 \qquad \hbox{and} \qquad \langle z^n, M_z^* 1\rangle_{\mathcal{H}_p(\mathbb{C})}=0.$$
\end{proof}

\begin{proposition}
\label{action5}
Let $p \in \mathbb{N}$. The adjoint operator of the integration operator in $ \mathcal{H}_p(\mathbb{C})$ is given by
\begin{equation}
\label{app}
I^* z^n= \frac{n^{2p}}{(n+1)^{2p}} z^{n-1}, \qquad I^*(1)=0.
\end{equation}
\end{proposition}
\begin{proof}
We have to prove that $ \langle I z^n, z^{m} \rangle_{\mathcal{H}_p(\mathbb{C})}=\langle  z^n, I^*z^{m} \rangle_{\mathcal{H}_p(\mathbb{C})}$, for $m$, $n \in \mathbb{N}_0$. We show first the result for $m \geq 1$, by \eqref{inner} we have 
$$ \langle  Iz^n, z^{m} \rangle_{\mathcal{H}_p(\mathbb{C})}= \frac{1}{n+1} \langle  z^{n+1}, z^{m} \rangle_{\mathcal{H}_p(\mathbb{C})}= \frac{1}{n+1} \frac{(n+1)!}{(n+2)^{2p}} \delta_{n+1, m}= \frac{n!}{(n+2)^{2p}} \delta_{n+1, m},$$		
whereas
$$ 	 \langle  z^n, I^* z^{m} \rangle_{\mathcal{H}_p(\mathbb{C})}= \frac{m^{2p}}{(m+1)^{2p}} \langle z^n, z^{m-1} \rangle_{\mathcal{H}_p(\mathbb{C})}= \frac{n!}{(n+2)^{2p}} \delta_{n,m-1}.$$
This shows the result for $m \geq 1$. For $m=0$ the result is trivial, indeed we have
$$
 \langle I 1, z^{m} \rangle_{\mathcal{H}_p(\mathbb{C})}=0 \qquad \hbox{and} \qquad
 \langle  1, I^*z^{m} \rangle_{\mathcal{H}_p(\mathbb{C})}= \frac{m^{2p}}{(m+1)^{2p}}\langle  1, z^{m-1} \rangle_{\mathcal{H}_p(\mathbb{C})}=0.$$
\end{proof}

We summarize the actions of all the adjoint operators previously studied in the following table

\begin{center}
	\begin{tabular}{| l | l |l |l|l|}
		\hline
		\rule[-4mm]{0mm}{1cm}
		{\bf Operators} & {$R_0^{*}z^n$} & {$\partial^* z^n$ } & {$M_z^{*}z^n$} &{$I^{*}z^n$}\\
		\hline
		\rule[-4mm]{0mm}{1cm}
		{} & $ \frac{(n+2)^{2p}}{(n+1)^{2p+1}} z^{n+1}$& $\frac{(n+2)^{2p}}{(n+1)^{2p}} z^{n+1}$ & $\frac{n^{2p+1}}{(n+1)^{2p}} z^{n-1}$& $\frac{n^{2p}}{(n+1)^{2p}} z^{n-1}$ \\
		\hline
		\hline
	\end{tabular}
\end{center}

In the following results we show that the previous adjoint operators can be written in peculiar forms. Precisely, let $\mathbf{A}^*$ be an adjoint operator in $ \mathcal{H}_p(\mathbb{C})$. We are able to show that the operators $R_0^*$, $\partial^*$, $M_z^*$ and $I^*$ can be written as the composition of the respective adjoint operators in the Fock space and a specific diagonal operator.

\begin{theorem}
\label{diag0}
Let $p \in \mathbb{N}$. Then the adjoint operator of the backward-shift operator can be written as
$$ R_0^*= I(\mathcal{I}+R_0I)^{2p},$$
where $ \mathcal{I}$ is the identity operator.
\end{theorem}
\begin{proof}
We show first the result for monomials $z^n$, with $n \in \mathbb{N}_0$.
By formula \eqref{action1} and the binomial theorem we have

\begin{eqnarray}
\nonumber
	R_0^{*}z^n&=&\frac{(n+2)^{2p}}{(n+1)^{2p+1}} z^{n+1}=\frac{1}{n+1} \left(1+ \frac{1}{n+1}\right)^{2p} z^{n+1}\\
			\label{a0}
	&=&\left( \sum_{k=0}^{2p} \binom{2p}{k} \frac{1}{(n+1)^{k+1}}\right) z^{n+1}.
\end{eqnarray}
Now we observe that for $k \geq 0$ we have
\begin{equation}
	\label{a2}
(R_0I)^k z^n= \frac{1}{(n+1)^k} z^{n}.
\end{equation}
This implies that
\begin{equation}
\label{a3}
I(R_0I)^{k}z^n= \frac{1}{(n+1)^{k+1}}z^{n+1}.
\end{equation}
By inserting \eqref{a3} in \eqref{a0} and applying another time the binomial theorem formula we get
$$ R_0^{*}z^n=\left( \sum_{k=0}^{2p} \binom{2p}{k} I(R_0I)^k\right) z^{n}=I\left( \sum_{k=0}^{2p} \binom{2p}{k} (R_0I)^k\right) z^{n}=I( \mathcal{I}+R_0I)^{2p} z^n.$$
Now, let us consider $f(z)= \sum_{n=0}^\infty a_n z^n$ and $g(z)= \sum_{n=0}^\infty b_n z^n$ in $ \mathcal{H}_p(\mathbb{C})$, with $ (a_n)_{n \in \mathbb{N}_0}$, $ (b_n)_{n \in \mathbb{N}_0} \subseteq \mathbb{C}$. By \eqref{a2} and the binomial theorem we get
\begingroup\allowdisplaybreaks	
\begin{eqnarray*}
\langle f,  I( \mathcal{I}+R_0I)^{2p} g\rangle_{\mathcal{H}_p(\mathbb{C})}&=&\left \langle f,  I( \mathcal{I}+R_0I)^{2p} \sum_{n=0}^\infty b_n z^n\right \rangle_{\mathcal{H}_p(\mathbb{C})}\\
&=& \left \langle f,   \sum_{n=0}^\infty \frac{(n+2)^{2p}}{(n+1)^{2p+1}} z^{n+1}b_n \right \rangle_{\mathcal{H}_p(\mathbb{C})}\\
&=& \left \langle f,   \sum_{n=1}^\infty  \frac{(n+1)^{2p}}{n^{2p+1}}z^{n}b_{n-1} \right \rangle_{\mathcal{H}_p(\mathbb{C})}\\
&=&  \sum_{n=1}^\infty  a_{n} \overline{b_{n-1}} \frac{n!}{(n+1)^{2p}}  \frac{(n+1)^{2p}}{n^{2p+1}}\\
&=&  \sum_{n=0}^\infty  a_{n+1} \overline{b_{n}} \frac{n!}{(n+1)^{2p}} \\
&=&  \left \langle \sum_{n=0}^\infty a_{n+1}z^n, g \right \rangle_{\mathcal{H}_p(\mathbb{C})} \\
&=& \langle R_0 f, g \rangle_{\mathcal{H}_p(\mathbb{C})}.
\end{eqnarray*}
\endgroup
This shows that $R_0^*= I(\mathcal{I}+R_0I)^{2p}$.

\end{proof}

\begin{theorem}
\label{diag1}
Let $p \in \mathbb{N}$. Then we can write the adjoint operators of the annihilation  operator in the following way
$$ \partial^*=M_z(\mathcal{I}+R_0I)^{2p},$$
where $ \mathcal{I}$ is the identity operator.
\end{theorem}
\begin{proof}
We show first the result for monomials $z^n$, with $n \in \mathbb{N}_0$.
By formula \eqref{action0} and the binomial theorem we have
\begin{equation}
	\label{ope0}
	\partial^* z^n= \left(\frac{n+2}{n+1}\right)^{2p} z^{n+1}= \left(1+ \frac{1}{n+1}\right)^{2p} z^{n+1}=\left( \sum_{k=0}^{2p} \binom{2p}{k} \frac{1}{(n+1)^k}\right) z^{n+1}.
\end{equation}
By \eqref{a2} we have that
\begin{equation}
\label{ope1}
M_z(R_0I)^{k} z^n= \frac{1}{(n+1)^k}z^{n+1}, \qquad k \geq 0.
\end{equation}
By \eqref{ope0}, \eqref{ope1} and the binomial theorem we obtain
$$ \partial^* z^n= M_z\left( \sum_{k=0}^{2p} \binom{2p}{k}(R_0I)^{k} \right) z^{n}=M_z (\mathcal{I}+IR_0)^{2p}.$$
Let us consider $f(z)= \sum_{n=0}^\infty a_n z^n$ and $g(z)= \sum_{n=0}^\infty b_n z^n$ in $ \mathcal{H}_p(\mathbb{C})$, with $ (a_n)_{n \in \mathbb{N}_0}$, $ (b_n)_{n \in \mathbb{N}_0} \subseteq \mathbb{C}$. By \eqref{a2} and the binomial theorem we get
\begin{eqnarray*}
	\langle f,  M_z( \mathcal{I}+R_0I)^{2p} g\rangle_{\mathcal{H}_p(\mathbb{C})}&=&\left \langle f,  M_z( \mathcal{I}+R_0I)^{2p} \sum_{n=0}^\infty b_n z^n\right \rangle_{\mathcal{H}_p(\mathbb{C})}\\
	&=& \left \langle f,   \sum_{n=0}^\infty  \left(\frac{n+2}{n+1}\right)^{2p} z^{n+1}b_n \right \rangle_{\mathcal{H}_p(\mathbb{C})}\\
	&=& \left \langle f,   \sum_{n=1}^\infty  \frac{(n+1)^{2p}}{n^{2p}} z^{n}b_{n-1} \right \rangle_{\mathcal{H}_p(\mathbb{C})}\\
	&=&  \sum_{n=1}^\infty  a_{n} \overline{b_{n-1}} \frac{n!}{n^{2p}} \\
	&=&  \left \langle \sum_{n=0}^\infty (n+1)a_{n+1}z^n, g \right \rangle_{\mathcal{H}_p(\mathbb{C})} \\
	&=& \langle \partial f, g \rangle_{\mathcal{H}_p(\mathbb{C})}.
\end{eqnarray*}
This shows that $\partial^* = M_z (\mathcal{I}+R_0I)^{2p}$.

\end{proof}

\begin{theorem}
\label{diag2}
Let $p \in \mathbb{N}$, then we can write the adjoint operator of the momentum operator as
\begin{equation}
	\label{ope}
	M_z^*= \partial(\mathcal{I}-R_0I)^{2p},
\end{equation}
where $ \mathcal{I}$ is the identity operator.
\end{theorem}
\begin{proof}
We show first the result for monomials $z^n$, with $n \in \mathbb{N}$. By the binomial theorem and Proposition \ref{action4} we have
\begin{eqnarray}
\nonumber
M_z^* z^n&=&\frac{n^{2p+1}}{(n+1)^{2p}} z^{n-1}=n \left(1- \frac{1}{n+1}\right)^{2p} z^{n-1}\\
		\label{bin}
&=&\left( \sum_{k=0}^{2p} \binom{2p}{k} (-1)^k \frac{n}{(n+1)^k}\right) z^{n-1}.
\end{eqnarray}
By \eqref{a2} we have that
\begin{equation}
	\label{ind}
	\partial  (R_0I)^k z^n= \frac{n}{(n+1)^k} z^{n-1}, \qquad k \geq 0.
\end{equation}
By inserting \eqref{ind} in \eqref{bin} we obtain
$$M_z^* z^n=\partial\left( \sum_{k=0}^{2p} \binom{2p}{k} (-R_0I)^k \right) z^{n}= \partial (\mathcal{I}-R_0I)^{2p}.$$
The result follows trivially for $n=0$.
Let $f(z)= \sum_{n=0}^\infty a_n z^n$ and $g(z)= \sum_{n=0}^\infty b_n z^n$ be functions in $ \mathcal{H}_p(\mathbb{C})$, with $ (a_n)_{n \in \mathbb{N}_0}$, $ (b_n)_{n \in \mathbb{N}_0} \subseteq \mathbb{C}$. By \eqref{a2} and the binomial theorem we get
\begingroup\allowdisplaybreaks	
\begin{eqnarray*}
	\langle f,  \partial (\mathcal{I}-R_0I)^{2p} g\rangle_{\mathcal{H}_p(\mathbb{C})}&=&\left \langle f,  \partial (\mathcal{I}-R_0I)^{2p} \sum_{n=0}^\infty b_n z^n\right \rangle_{\mathcal{H}_p(\mathbb{C})}\\
	&=& \left \langle f,   \sum_{n=1}^\infty  \frac{n^{2p+1}}{(n+1)^{2p}} z^{n-1}b_n \right \rangle_{\mathcal{H}_p(\mathbb{C})}\\
	&=& \left \langle f,   \sum_{n=0}^\infty  \frac{(n+1)^{2p+1}}{(n+2)^{2p}} z^{n}b_{n+1} \right \rangle_{\mathcal{H}_p(\mathbb{C})}\\
	&=&  \sum_{n=0}^\infty  a_{n} \overline{b_{n+1}} \frac{(n+1)!}{(n+2)^{2p}} \\
	&=&  \left \langle \sum_{n=1}^\infty a_{n-1}z^n, g \right \rangle_{\mathcal{H}_p(\mathbb{C})} \\
	&=& \langle M_z f, g \rangle_{\mathcal{H}_p(\mathbb{C})}.
\end{eqnarray*}
\endgroup
This shows that $M_z^* = \partial  (\mathcal{I}-R_0I)^{2p}$.

\end{proof}

\begin{theorem}
	\label{diag3}
Let $p \in \mathbb{N}$. Then we can write the adjoint of the integration operator as
\begin{equation}
	\label{app2}
	I^*=R_0(\mathcal{I}-R_0I)^{2p},
\end{equation}
where $ \mathcal{I}$ is the identity operator.
\end{theorem}
\begin{proof}
We show first the result for monomials $z^n$, with $n \in \mathbb{N}$. By using  the binomial theorem and formula \eqref{app} we get
\begin{equation}
	\label{int}
	I^{*} z^n=\frac{n^{2p}}{(n+1)^{2p}}z^{n-1}= \left(1- \frac{1}{n+1}\right)^{2p}z^{n-1}=\left(\sum_{k=0}^{2p} \binom{2p}{k} (-1)^k \frac{1}{(n+1)^k} \right) z^{n-1}.
\end{equation}
By \eqref{a2} we have that
\begin{equation}
	\label{ind2}
	R_{0}(R_0I)^{k} z^{n}= \frac{1}{(n+1)^{k}} z^{n-1}, \qquad k \geq 0.
\end{equation}

By inserting \eqref{ind2} in \eqref{int} and using again the binomial theorem we get
$$ I^{*} z^n= R_{0}\left(\sum_{k=0}^{2p} \binom{2p}{k} (-R_0I)^{k}  \right) z^{n}=R_{0}(\mathcal{I}-R_0I)^{2p}.$$
The result follows trivially for $n=0$.
Let $f(z)= \sum_{n=0}^\infty a_n z^n$ and $g(z)= \sum_{n=0}^\infty b_n z^n$ be functions in $ \mathcal{H}_p(\mathbb{C})$, with $ (a_n)_{n \in \mathbb{N}_0}$, $ (b_n)_{n \in \mathbb{N}_0} \subseteq \mathbb{C}$. By the binomial theorem we get
\begingroup\allowdisplaybreaks	
\begin{eqnarray*}
	\langle f,  R_0 (\mathcal{I}-R_0I)^{2p} g\rangle_{\mathcal{H}_p(\mathbb{C})}&=&\left \langle f,  R_0 (\mathcal{I}-R_0I)^{2p} \sum_{n=0}^\infty b_n z^n\right \rangle_{\mathcal{H}_p(\mathbb{C})}\\
	&=& \left \langle f,   \sum_{n=1}^\infty  \frac{n^{2p}}{(n+1)^{2p}} z^{n-1}b_n \right \rangle_{\mathcal{H}_p(\mathbb{C})}\\
	&=& \left \langle f,   \sum_{n=0}^\infty  \frac{(n+1)^{2p}}{(n+2)^{2p}} z^{n}b_{n+1} \right \rangle_{\mathcal{H}_p(\mathbb{C})}\\
	&=&  \sum_{n=0}^\infty  a_{n} \overline{b_{n+1}} \frac{n!}{(n+2)^{2p}} \\
	&=&  \left \langle \sum_{n=1}^\infty \frac{a_{n-1}}{n}z^n, g \right \rangle_{\mathcal{H}_p(\mathbb{C})} \\
	&=& \langle I f, g \rangle_{\mathcal{H}_p}.
\end{eqnarray*}
\endgroup
This shows that $I^{*} = R_{0}(\mathcal{I}-R_0I)^{2p}$.

\end{proof}

We sum up all the previous results in the following table

\begin{center}
	\begin{tabular}{| l | l |l |l|l|}
		\hline
		\rule[-4mm]{0mm}{1cm}
		{\bf Operators} & {$R_0^{*}$} & {$\partial^*$ } & {$M_z^{*}$} &{$I^{*}$}\\
		\hline
		\rule[-4mm]{0mm}{1cm}
		{} & $ I(\mathcal{I}+R_0I)^{2p}$& $M_z (\mathcal{I}+R_0I)^{2p}$ & $ \partial  (\mathcal{I}-R_0I)^{2p}$& $ R_{0}(\mathcal{I}-R_0I)^{2p}$ \\
		\hline
		\hline
	\end{tabular}
\end{center}

\begin{remark}
If we consider $p=0$ in Theorem \ref{diag0}, Theorem \ref{diag1}, Theorem \ref{diag2} and Theorem \ref{diag3} we get the results that hold in classical theory of the Fock spaces.
\end{remark}

Unlike to what happens in the classic case the operators $M_z$ and $M_{z}^*$ do not satisfy a commutation relation, as well as $R_0$ and $R_0^*$, see \eqref{comm} and \eqref{com}.
We show that the commutators in the space $ \mathcal{H} _p(\mathbb{C})$ are related to peculiar diagonal operators. We denote a generic diagonal operator by
$$ D= \hbox{diag} (d_0, d_1, d_2,...).$$
We define the forward and backwards shifts of $D$ as
\begin{equation}
	\label{forw}
	D^{(1)}=\hbox{diag}(0, d_0, d_1,...),
\end{equation}
\begin{equation}
\label{back}
D^{(-1)}=\hbox{diag}(d_1, d_2,d_3,...),
\end{equation}
respectively.
\begin{proposition}
\label{comm2}
The operators $M_z$ and $M_z^{*}$ satisfy the following  formal relation
\begin{equation}
\label{comm8}
[M_z, M_z^*]=\left[\mathsf{A} \left(Id- [\mathsf{A}^{(-1)}]^{-1}\right)^{2p}-\mathsf{A}^{(-1)}\left(Id- [\mathsf{A}^{(-2)}]^{-1}\right)^{2p}\right],
\end{equation}
where $Id$ is the identity matrix, $ \mathsf{A}:= \hbox{diag}(0,1,2,3,...)$ and $[.]^{-1}$ is the inverse matrix.
\end{proposition}
\begin{proof}
We prove first the result on the monomials. By Theorem \ref{diag2} we have
\begin{eqnarray*}
[M_z, M_z^*](z^n)&=&M_z  M_z^*(z^n)-M_{z}^* M_z (z^n) \\
&=& \left(\sum_{k=0}^{2p} \binom{2p}{k}(-1)^k M_z \partial (R_0I)^k- \sum_{k=0}^{2p} \binom{2p}{k}(-1)^k \partial (R_0I)^k M_z\right)z^n.
\end{eqnarray*}
By \eqref{a2} we have
$$ M_z \partial (IR_0)^k(z^n)= \frac{n}{(n+1)^k}z^n \qquad \hbox{and} \qquad \partial (IR_0)^k M_z (z^n)= \frac{(n+1)}{(n+2)^k}z^n.$$
Then by the binomial formula we have
\begin{eqnarray}
	\nonumber
	[M_z, M_z^*](z^n)&=&\left(n\sum_{k=0}^{2p} \binom{2p}{k}(-1)^k  \frac{1}{(n+1)^k}-(n+1)\sum_{k=0}^{2p} \binom{2p}{k}(-1)^k  \frac{1}{(n+2)^k}\right)z^n\\
\label{comm7}
	&=& \left[n \left(1- \frac{1}{n+1}\right)^{2p}-(n+1)\left(1- \frac{1}{n+2}\right)^{2p}\right]z^n.
\end{eqnarray}
Now, by setting
\begin{equation}
\label{not}
 \mathsf{A}:=\hbox{diag}(0,1,2,3,...,n,...), \qquad \mathsf{B}:=(1,2,3,...,n+1,...)
\end{equation}
\begin{equation}
\label{not1}
\mathsf{C}:=\hbox{diag}\left(1, \frac{1}{2}, \frac{1}{3},..., \frac{1}{n+1},...\right) \qquad \mathsf{D}:=\hbox{diag}\left(\frac{1}{2}, \frac{1}{3}, \frac{1}{4},..., \frac{1}{n+2},...\right)
\end{equation}
we can formally write \eqref{comm7} in the following way
$$ [M_z, M_z^*]=\left[\mathsf{A} \left(Id- \mathsf{C}\right)^{2p}-\mathsf{B}\left(Id- \mathsf{D}\right)^{2p}\right].$$
Finally we observe that $ \mathsf{B}= \mathsf{A}^{(-1)}$, $ \mathsf{C}=[\mathsf{A}^{(-1)}]^{-1}$ and $ \mathsf{D}=[\mathsf{A}^{(-2)}]^{-1}$. These considerations imply \eqref{comm8}.

\end{proof}
\begin{remark}
	If we consider $p=0$ in Proposition \ref{comm2} we have that $M_z^{*}=\partial$ and
	$$ \mathsf{A}-\mathsf{A}^{(-1)}= \hbox{diag}(-1,-1,-1,...)=-id.$$
	Hence we get back to formula \eqref{comm}.
\end{remark}
\begin{proposition}
\label{comm3}
The operators $R_0$ and $R_0^*$ satisfy the following formal relation
\begin{equation}
\label{comm9}
[R_0, R_0^*]=D_0 \left[\mathsf{A}^{(-1)} \left(Id+\mathsf{E}\right)^{2p}-\mathsf{A} \left(Id+ \mathsf{E}^{(-1)}\right)^{2p}\right],
\end{equation}
where $Id$ is the identity matrix, $ \mathsf{A}:= \hbox{diag}(0,1,2,3,...)$ , $\mathsf{E}:= \hbox{diag}\left(1,1, \frac{1}{2}, \frac{1}{3},..., \frac{1}{n},...\right)$ and $D_0$ is the diagonal matrix defined in \eqref{mat}.
\end{proposition}
\begin{proof}
We prove the result on the monomials. By Theorem \ref{diag0} and for $n \geq 1$ we have
\begin{eqnarray}
\nonumber
[R_0, R_0^*](z^n)&=& R_0 R_0^*(z^n)-R_0^*R_0(z^n)\\
\label{comm1}
&=& \sum_{k=0}^{2p} \binom{2p}{k} (R_0I)^{k+1}(z^n)- \sum_{k=0}^{2p} \binom{2p}{k} I(R_{0}I)^k R_0(z^n).
\end{eqnarray}
By formula \eqref{a2} we have
$$ (R_0I)^{k+1}(z^n)= \frac{z^n}{(n+1)^{k+1}}, \qquad \hbox{and} \qquad [I(R_{0}I)^k R_0](z^n)= \frac{z^n}{n^{k+1}}.$$
Hence we get
\begin{eqnarray*}
[R_0, R_0^*](z^n)&=&\left( \sum_{k=0}^{2p} \binom{2p}{k} \frac{1}{(n+1)^{k+1}}- \sum_{k=0}^{2p} \binom{2p}{k} \frac{1}{n^{k+1}} \right) z^n\\
&=&- \frac{1}{n(n+1)} \left(\sum_{k=0}^{2p}\binom{2p}{k} \frac{(n+1)^{k+1}-n^{k+1}}{n^k(n+1)^k}\right)\\
&=& - \frac{1}{n(n+1)}\left( (n+1)\sum_{k=0}^{2p}\binom{2p}{k} \frac{1}{n^k}-n \sum_{k=0}^{2p}\binom{2p}{k} \frac{1}{(n+1)^k} \right)z^n.
\end{eqnarray*}
From the binomial theorem we get
$$ [R_0, R_0^*](z^n)=- \frac{1}{n(n+1)} \left((n+1) \left(1+ \frac{1}{n}\right)^{2p}-n \left(1+ \frac{1}{n+1}\right)^{2p}\right)z^n$$
Now, for $n=0$ we have that $R_0^*(1)=4^p$ therefore by \eqref{comm1} we obtain
$$ [R_0, R_0^*](1)= R_0 R_0^*(1)-R_0^*R_0(1)=4^p.$$
Therefore we proved that
\begin{equation}
\label{res}
[R_0, R_0^*](z^n)= \begin{cases}
	4^p, \quad n=0\\
	- \frac{1}{n(n+1)}\left((n+1) \left(1+ \frac{1}{n}\right)^{2p}-n \left(1+ \frac{1}{n+1}\right)^{2p}\right)z^n, \quad n=1,2,3,... 
\end{cases}
\end{equation}
Now by setting
$$ \mathsf{E}:= \hbox{diag}\left(1,1, \frac{1}{2}, \frac{1}{3},..., \frac{1}{n},...\right)$$
and by the notations introduced in \eqref{not} and \eqref{not1} we can write \eqref{res} formally as 
$$ [R_0, R_0^*]= D_0 \left[\mathsf{B}\left(Id+\mathsf{E}\right)^{2p}-\mathsf{A} \left(Id+ \mathsf{C}\right)^{2p}\right].$$
Finally, by observing that $ \mathsf{B}= \mathsf{A}^{(-1)}$ and $ \mathsf{C}=\mathsf{E}^{(-1)}$ we get \eqref{comm9}.
\end{proof}

\begin{remark}
If we take $p=0$ in Proposition \ref{comm3} we have that $R_0^{*}=I$ and
$$ \mathsf{A}^{(-1)}- \mathsf{A}=\hbox{diag}(1,1,1,...)=Id.$$
Since $D_0 Id=D_0$ we get back the identity \eqref{comm5}.
\end{remark}

\subsection{The $\mathcal{H}_p$-Bargmann transform}

In this section we show another characterization of the space $ \mathcal{H}_p(\mathbb{C})$. This is obtained by means of a peculiar map from $L^2(\mathbb{R})$ onto $ \mathcal{H}_p(\mathbb{C})$.
\begin{proposition}
\label{kernel1}
Let $p \in \mathbb{N}$. For every $x \in \mathbb{R}$ and $z \in \mathbb{C}$ we define the function
\begin{equation}
A_{p}(z,x)= \sum_{n=0}^{\infty} \frac{(n+1)^{p} z^n}{\sqrt{n!}} \xi_n(x),
\end{equation}
where $\xi_n$ are the normalized Hermite  functions. Then we have:
\begin{itemize}
\item[1)] the function $A_p(.,x)$ is entire for every $x \in \mathbb{R}$.
\item[2)] a function $f \in \mathcal{H}_p(\mathbb{C})$ if and only if there exists $\varphi \in L^{2}(\mathbb{R})$ such that
\begin{equation}
f(z)= \int_{\mathbb{R}} A_p(z,x) \varphi(x) dx= \langle \varphi, \overline{A_p(z,.)} \rangle_{L^2}.
\end{equation}
\end{itemize}
\end{proposition}
\begin{proof}
We know by \cite{Hille} that the normalized Hermite functions are uniformly bounded by some constants, i.e.
$$ \exists C>0 \quad \hbox{such that} \qquad | \xi_n(x)|<C \quad \hbox{for every} \quad n \in \mathbb{N} \quad \hbox{and} \quad x \in \mathbb{R}.$$
This implies that the series in \eqref{kernel1} is convergent, and so $A_p(., x)$ is entire. 
\\Now, we show the second point of the statement. Let us suppose that $f(z)=\langle g, \overline{A_p(z,.)} \rangle_{L^2}$, for $\varphi \in L^2(\mathbb{R})$. We have to prove that $f \in \mathcal{H}_p(\mathbb{C})$. Then we have
\begin{eqnarray*}
 f(z)&=&\int_{\mathbb{R}} \left( \sum_{n=0}^\infty \frac{(n+1)^p z^n}{\sqrt{n!}} \xi_n(x) \varphi(x)\right) dx\\
 &=& \sum_{n=0}^{\infty}\frac{(n+1)^p z^n}{\sqrt{n!}} \int_{\mathbb{R}} \xi_n(x) \varphi(x)dx\\
 &=&  \sum_{n=0}^{\infty} z^n \alpha_n,
\end{eqnarray*}
where $ \alpha_n:= \frac{(n+1)^p}{\sqrt{n!}} \int_{\mathbb{R}} \xi_n(x) \varphi(x)dx$.
Since $ \{\xi_n\}_{n=0}^\infty$ is an orthonormal basis of $L^2(\mathbb{R})$ and by the Parseval's equality we have that
\begin{eqnarray*}
\sum_{n=0}^\infty  \frac{n!}{(n+1)^{2p}} | \alpha_n|^2&=& \frac{1}{\sqrt{2 \pi}} \sum_{n=0}^\infty \left| \int_{\mathbb{R}} \xi_n(x) \varphi(x) dx \right|^2\\
&=& \int_{\mathbb{R}} |\varphi(x)|^2 dx\\
&=&\| \varphi\|_{L^2}^2< \infty.
\end{eqnarray*}
This proves that $f \in \mathcal{H}_p(\mathbb{C})$. 
\\Now, we show the converse. From the definition of the space $ \mathcal{H}_p(\mathbb{C})$ we have $f(z)= \sum_{n=0}^\infty z^n a_n$ and $ \sum_{n=0}^\infty\frac{n!}{(n+1)^{2p}}|a_n|^2< \infty$, with $ (a_n)_{n \in \mathbb{N}_0} \subseteq \mathbb{C}$. We set
$$ \varphi(x):= \sum_{n=0}^\infty \frac{\sqrt{n!}}{(n+1)^{p}} a_n \xi_n(x).$$
This function belongs to $L^2(\mathbb{R})$. Indeed we have
$$ \| \varphi\|_{L^2}^2= \sum_{n=0}^\infty \frac{n!}{(n+1)^{2p}}|a_n|^2< \infty.$$
Since the normalized Hermite functions are a basis of the space $L^2(\mathbb{R})$ we have
$$ \langle \varphi, \overline{A_p(z,.)} \rangle_{L^2}=\sum_{n=0}^\infty  \frac{z^n \sqrt{n!}}{(n+1)^{p}} \frac{(n+1)^{p}}{\sqrt{n!}} a_n= \sum_{n=0}^{\infty} z^n a_n=f(z).$$
This concludes the proof.
\end{proof}

The previous result gives the motivation to consider a Bargmann transform defined in the space $ \mathcal{H}_p(\mathbb{C})$.
\begin{definition}
\label{Fbarg}
Let $p \in \mathbb{N}$. For any $\varphi \in L^2(\mathbb{R})$ we define the $\mathcal{H}_p$-Bargmann transform as
$$\mathcal{B}_p(\varphi(z)):= \int_{\mathbb{R}}  A_p(z,x) \varphi(x)dx, \qquad A_p(z,x):= \frac{1}{(2\pi)^{\frac{1}{4}}} \sum_{n=0}^\infty \frac{(n+1)^p z^n}{n!} h_n(x),$$
where $h_n$ are the Hermite functions defined in \eqref{harmite}.
\end{definition}

\begin{remark}
If we consider $p=0$ in Definition \ref{Fbarg} we get back  to the classical Bargmann transform, see \eqref{kernel} and \eqref{genkern}.
\end{remark}
To show the next result we need the following notation
\begin{equation}
\label{nota1}
 A_{p}^z:=A_p(.,x),
\end{equation}
that it means that the variable $z$ is fixed.
\begin{proposition}
Let $p \in \mathbb{N}$. For any $z$, $w \in \mathbb{C}$ we have that
\begin{equation}
\label{nota}
\langle A_{p}^z, A_p^w \rangle_{L^2}= \frac{e^{z \bar{w}}}{z \bar{w}} T_{2p+1}(z \bar{w}).
\end{equation}
\end{proposition}
\begin{proof}
Let $z$, $w \in \mathbb{C}$, by the definition of the kernel $A_p$ we have
\begin{eqnarray}
\nonumber
\langle A_{p}^z, A_p^w \rangle_{L^2}&=& \int_{\mathbb{R}} A_p^z(x) \overline{A_p^w(x)} dx\\
\label{ref1}
&=& \int_{\mathbb{R}} \left( \sum_{j=0}^{\infty} \frac{(j+1)^p z^j}{\sqrt{j!}} \xi_j(x)\right)\left( \sum_{k=0}^{\infty} \frac{(k+1)^p \bar{w}^k}{\sqrt{k!}} \xi_k(x)\right) dx.
\end{eqnarray}
The Hermite functions form an orthogonal basis of the space $L^2(\mathbb{R})$, hence we have
\begin{equation}
\label{int2}
\int_{\mathbb{R}} \xi_j(x) \xi_k(x)  dx=\delta_{j,k}, \qquad j,k \in \mathbb{N},
\end{equation}
where $\delta_{k,j}$ is the Kronecker symbol. By Proposition \ref{kernel1} we can exchange the double series in \eqref{ref1} and so by formula \eqref{int2} we have
\begin{eqnarray*}
\langle A_{p}^z, A_p^w \rangle_{L^2}&=& \sum_{j,k=0}^{\infty} \frac{(j+1)^p (k+1)^p z^j \bar{w}^k}{ \sqrt{j! k!}} \int_{\mathbb{R}} \xi_{j}(x) \xi_k(x) dx\\
&=& \sum_{k=0}^{\infty} \frac{(k+1)^{2p} z^k \bar{w}^k}{ k!}\\
&=& K_p(z,w).
\end{eqnarray*}
The result follows by Proposition \ref{rel}.
\end{proof}
\begin{remark}
Formula \eqref{nota} gives another representation of the Touchard polynomials in terms of the $L^2$-inner product of the kernel of the $\mathcal{H}_p$-Bargmann transform.
\end{remark}

Now, our goal is to find a closed expression of of the $\mathcal{H}_p$-Bargmann transform. We recall that the number operator is defined by means of the combination of the momentum and annihilation operators, i.e.  $ M_z \partial$.
It acts on the monomial $z^n$, with $ n\geq0$, in the following way
$$ \left(M_z \partial\right)z^n=nz^n.$$
By iterating the action of the number operator to $z^n$ we get
\begin{equation}
\label{closed1}
\left(M_z \partial\right)^jz^n=n^jz^n, \qquad j \geq 0.
\end{equation}
In \cite{KNO} the author proved that the action of the number operator can be written in terms of the Sterling numbers of the second type, see \eqref{Stirling}. Precisely it holds that
\begin{equation}
\label{closed2}
\left(M_z \partial\right)^j= \sum_{k=0}^j S(j,k) z^k \frac{d^k}{d z^k}.
\end{equation}
\begin{theorem}
\label{mainBarg}
For $p \in \mathbb{N}$ we can write the kernel of the $\mathcal{H}_p$-Bargmann transform as
\begin{equation}
A_p(z,x)= \left(\mathcal{I}+M_z \partial\right)^p A(z,x), \qquad \forall (z,x) \in \mathbb{C} \times \mathbb{R}.
\end{equation}
Moreover, we can write
\begin{equation}
	\label{secine}
	A_{p}(z,x)= \left( \partial M_z\right)^p A(z,x), 
\end{equation}
where $A(z,x)$ is the kernel of the classical Bargmann transform, see \eqref{kernel}.
\end{theorem}
\begin{proof}
By the binomial formula we have
$$ (n+1)^p= \sum_{j=0}^p \binom{p}{j} n^j.$$
This implies that we can write the kernel of the $\mathcal{H}_p$-Bargmann transform as
\begin{eqnarray}
\nonumber
A_p(z,x)&=& \frac{1}{(2\pi)^{\frac{1}{4}}} \sum_{n=0}^{\infty} \frac{(n+1)^p}{n!} h_n(x) z^n\\
\nonumber
&=&\frac{1}{(2\pi)^{\frac{1}{4}}}\sum_{n=0}^{\infty} \frac{1}{n!} \left( \sum_{j=0}^{p} \binom{p}{j} n^j\right) h_n(x) z^n\\
\label{closed3}
&=& \sum_{j=0}^p \binom{p}{j} \left[\frac{1}{(2\pi)^{\frac{1}{4}}} \sum_{n=0}^\infty \frac{h_n(x)}{n! } n^j z^n\right].
\end{eqnarray}
By \eqref{closed1} and using another time the binomial formula we have
\begin{eqnarray*}
A_{p}(z,x)&=& \sum_{j=0}^p \binom{p}{j} \left[ \frac{1}{(2\pi)^{\frac{1}{4}}} \sum_{n=0}^\infty \frac{h_n(x)}{n!} \left(M_z \partial\right)^j z^n  \right]\\
&=& \sum_{j=0}^p \binom{p}{j}\left(M_z \partial\right)^j \left[ \frac{1}{(2\pi)^{\frac{1}{4}}} \sum_{n=0}^\infty \frac{h_n(x)}{n! }  z^n  \right]\\
&=& \sum_{j=0}^p \binom{p}{j}\left(M_z \partial\right)^j A(z,x)\\
&=& \left(\mathcal{I}+M_z \partial\right)^p A(z,x).
\end{eqnarray*}
Finally the equality \eqref{secine} follows from \eqref{comm}.
\end{proof}

A direct application of the previous result is given by the following generating formula.
\begin{corollary}
For $p \in \mathbb{N}$ we have that
$$ \sum_{p=0}^{\infty} \frac{A_p(z,x)}{p!}=exp \left[\left(\mathcal{I}+M_z \partial\right) \right]A(z,x).$$
\end{corollary}

Now, we show a recurrence-like relation for the kernel $A_p(z,x)$.
\begin{proposition}
Let $p \in \mathbb{N}$. Then we have
$$ A_{p+1}(z,x)-A_p(z,x)=M_z \left(\mathcal{I}+ \partial M_z\right)^p \partial A(z,x), \qquad \forall (z,x) \in \mathbb{C} \times \mathbb{R}.$$
\end{proposition}
\begin{proof}
By Theorem \ref{mainBarg} we have
\begin{eqnarray*}
A_{p+1}(z,x)&=&\left( \mathcal{I}+ M_z \partial\right)^{p+1} A(z,x)\\
&=& \left( \mathcal{I}+ M_z \partial\right) A_p(z,x)\\
&=& A_p(z,x)+ M_z \partial A_p(z,x).
\end{eqnarray*}
By using another time Theorem \ref{mainBarg} we have
\begin{equation}
\label{rec}
A_{p+1}(z,x)-A_p(z,x)=M_z \partial \left(\mathcal{I}+M_z \partial\right)^p A(z,x).
\end{equation}
Now, for $j \geq 0$ we have
\begin{equation}
\label{ref2}
\partial\left(M_z \partial\right)^j=\partial\left(M_z \partial\right)^{j-1}M_z \partial=\left(\partial M_z\right)^j \partial.
\end{equation}

By the binomial formula and \eqref{ref2} we have
\begin{eqnarray}
\nonumber
\partial \left(\mathcal{I}+M_z \partial \right)^pA(z,x)&=& \sum_{j=0}^{p} \binom{p}{j} \partial\left(M_z \partial\right)^j A(z,x)\\
\nonumber
&=&\sum_{j=0}^{p} \binom{p}{j} \left(\partial M_z\right)^j \partial A(z,x)\\
\label{ref3}
&=& \left(\mathcal{I}+ \partial M_z\right)^p\partial A(z,x)
\end{eqnarray}
By plugging \eqref{ref3} into \eqref{rec} we get the final result.
\end{proof}

Now, we prove that the kernel of $\mathcal{B}_p$ can be written also in terms of the Sterling numbers of the second kind.
\begin{proposition}
Let $p \in \mathbb{N}$. Then we can write the kernel of the $\mathcal{H}_p$-Bargmann transform as
\begin{equation}
\label{closedB}
A_p(z,x)= \sum_{j=0}^p \sum_{k=0}^j \binom{p}{j} S(j,k) z^k  \frac{d^k}{dz^k}A(z,x),
\end{equation}
where $S(n,k)$ are the Sterling number of the second type.
\end{proposition}
\begin{proof}
By Theorem \ref{mainBarg} and formula \eqref{closed2} we have
$$A_p(z,x)=\sum_{j=0}^p \binom{p}{j}\left(z \partial\right)^j A(z,x)=\sum_{j=0}^p \sum_{k=0}^j \binom{p}{j} S(j,k) z^k  \frac{d^k}{dz^k}A(z,x).$$
This proves the result.
\end{proof}
To show the next result we recall the extension of the Hermite polynomials in the complex variable $z$, see \cite{S, V}. Precisely, these are given by 

$$ H_n(z)= \sum_{k=0}^{\left\lfloor \frac{n}{2} \right\rfloor }\frac{n! (-1)^k}{k! (n-2k)!} \frac{(\sqrt{2}z)^{n-2k}
}{2^{\frac{n}{2}}}.$$
These polynomials satisfied the following recurrence relation
\begin{equation}
\label{recurr}
H_{k+1}(z)=z H_k(z)-\frac{d}{dz} H_k(z).
\end{equation}

\begin{proposition}
	\label{der}
For $k \geq 0$ we have
$$
\frac{d^k}{dz^k}A(z,x)=H_k(x-z)A(z,x),
$$
where $A(z,x)$ is the kernel of the classical Bargmann transform.
\end{proposition}
\begin{proof}
We prove the result by induction on $k$. For $k=0$ the result is trivial because $H_0(x-z)=1$. Now we suppose that it is true for $k$ and we prove it for $k+1$. 
We start by observing that
$$ \frac{d}{dz}A(z,x)=(x-z)A(z,x).$$
Hence by \eqref{recurr} and the inductive hypothesis we get
\begin{eqnarray*}
\frac{d^{k+1}}{dz^{k+1}}A(z,x)&=& \frac{d}{dz} \left( \frac{d^k}{dz^k}A(z,x)\right)\\
&=& \frac{d}{dz} \left(H_k(x-z)A(z,x) \right)\\
&=& -\frac{d}{dz}H_k(x-z)A(z,x)+H_k(x-z)(x-z)A(z,x)\\
&=& \left[H_k(x-z)(x-z)-\frac{d}{dz} H_k(x-z)\right]A(z,x)\\
&=& H_{k+1}(x-z)A(z,x).
\end{eqnarray*}
\end{proof}

The above result paved the way to get a closed formula for the $\mathcal{H}_p$-Bargmann transform.

\begin{proposition}
\label{barp}
Let $p \in \mathbb{N}$. Then we can write the $\mathcal{H}_p$-Bargmann transform in terms of the classical Bargmann transform. Precisely for any $ \varphi \in L^2(\mathbb{R})$ it holds that
\begin{eqnarray}
\nonumber
\mathcal{B}_p(\varphi(z))&=& \left(\mathcal{I}+M_z \partial\right)^p (\mathcal{B}\varphi)(z)\\
\nonumber
&=& \left( \partial M_z\right)^p (\mathcal{B}\varphi)(z)\\
\label{derbar}
&=& \sum_{j=0}^p \sum_{k=0}^{j} \binom{p}{j}S(j,k) z^k \frac{d^k}{dz^k} (\mathcal{B}\varphi)(z).
\end{eqnarray}
Moreover, we can write the $\mathcal{H}_p$-Bargmann transform as 
\begin{equation}
\label{closed12}
\mathcal{B}_p(\varphi(z))= \left(\frac{1}{2 \pi}\right)^{\frac{1}{4}}\sum_{j=0}^{p} \sum_{k=0}^{j} \binom{p}{j}S(j,k)z^k \int_{\mathbb{R}} H_k(x-z) e^{- \frac{x^2}{4}+xz- \frac{z^2}{2}} \varphi(x) dx.
\end{equation}
\end{proposition}
\begin{proof}
	The result follows from Theorem \ref{mainBarg} and Proposition \ref{closedB}. Finally by \eqref{derbar} and Proposition \ref{der} we get \eqref{closed12}.
\end{proof}

\subsection{The Fock space connected to the $ \mathcal{H}_p$ space}

In this section we study a connection between the space $ \mathcal{H}_p$ and the Fock space. This is achieved by a  diagonal specific operator. A similar operator has been also considered in \cite[Prop. 2.6]{Ner}, but in this paper we follow a different approach. At the end of this section, by using the connection between the spaces $\mathcal{F}(\mathbb{C})$ and $\mathcal{H}_p(\mathbb{C})$, we show also some basic properties of the $\mathcal{H}_p$-Bargmann transform.

\begin{theorem}
\label{uni}
Let $p \in \mathbb{N}$. The the operator 
\begin{equation}
	\label{ope2}
	\mathcal{E}_p:=\left(\mathcal{I}+M_z \partial\right)^p.
\end{equation}
realizes a surjective isometry from the space $ \mathcal{F}(\mathbb{C})$ onto $ \mathcal{H}_p(\mathbb{C})$.
\end{theorem}
\begin{proof}
We show first that the operator $ \mathcal{E}_p$ is injective. To do this we prove that $ \mathcal{E}_p$ realizes an isometry. Let us consider $f \in \mathcal{F}(\mathbb{C})$. We apply the operator $ \mathcal{E}_p$ to a function $f=\sum_{n=0}^{\infty}z^n a_n$, with $(a_n)_{n \in \mathbb{N}_0}\subseteq \mathbb{C}$. Since $ \mathcal{E}_p(z^n)=(n+1)^p z^n$ we get
\begin{eqnarray}
	\nonumber
	\mathcal{E}_p(f)&=& \mathcal{E}_p \left( \sum_{n=0}^{\infty} z^n a_n\right)\\
	\nonumber
	&=& \sum_{n=0}^{\infty} (n+1)^p z^n a_n\\
	\label{app1}
	&=& \sum_{n=0}^{\infty} z^n \alpha_n,
\end{eqnarray}
where $\alpha_n:=(n+1)^p a_n$. Therefore we have
$$ \| \mathcal{E}_p(f) \|_{\mathcal{H}_p(\mathbb{C})}^2=\sum_{n=0}^\infty \frac{n!}{(n+1)^{2p}}| \alpha_n|^2=\sum_{n=0}^\infty n! | a_n|^2=\| f \|_{\mathcal{F}(\mathbb{C})}^2.$$
This proves that the operator $ \mathcal{E}_p$ is an isometry between the spaces  $ \mathcal{F}(\mathbb{C})$ and $ \mathcal{H}_p(\mathbb{C})$. Now, we show that the operator $ \mathcal{E}_p$ is surjective. Let $g \in \mathcal{H}_p(\mathbb{C})$. Hence we can write $g(z)=\sum_{n=0}^{\infty} z^na_n$, with $ (a_n)_{n \in \mathbb{N}_0} \subset \mathbb{C}$. To achieve our aim we have to find a function $f \in \mathcal{F}(\mathbb{C})$, such that
$$ g(z)= \mathcal{E}_p(f)(z).$$
Precisely we have to pick a sequence $(b_n)_{n \in \mathbb{N}_0} \subset \mathbb{C}$ such that $f(z)=\sum_{n=0}^{\infty} z^n b_n$ and $\| f\|_{\mathcal{F}}^2=\sum_{n=0}^{\infty}n! |b_n|^2< \infty.$ Let us consider $b_n=\frac{a_n}{(n+1)^p}$. Then we have
\begin{eqnarray*}
\mathcal{E}_p(f)(z)&=& \sum_{n=0}^{\infty}(n+1)^p z^n b_n\\
&=& \sum_{n=0}^{\infty} z^n a_n\\
&=& g(z).
\end{eqnarray*}
Moreover since $g \in \mathcal{H}_p(\mathbb{C})$ we have
$$
\| f\|_{\mathcal{F}(\mathbb{C})}^2= \sum_{n=0}^{\infty} n! |b_n|^2= \sum_{n=0}^\infty \frac{n!}{(n+1)^{2p}}|a_n|^2 < \infty.$$
Therefore, with the above choice of $(b_n)_{n \in \mathbb{N}_0}$ we have that $f \in \mathcal{F}(\mathbb{C})$.
\end{proof}

\begin{remark}
The operator $\mathcal{E}_p$ is the identity if we take $p=0$.
\end{remark}

\begin{proposition}
\label{adj3}
For $p\in \mathbb{N}$ the adjoint of the operator $ \mathcal{E}_p$ is given by $(R_0I)^p$. Moreover the inverse of the operator $ \mathcal{E}_p$ is given by its adjoint.
\end{proposition}
\begin{proof}
Let us consider $g \in \mathcal{H}_p( \mathbb{C})$ and $f \in \mathcal{F}(\mathbb{C})$. To prove the result we have to show that
$$ \langle \mathcal{E}_p(f), g \rangle_{\mathcal{H}_p(\mathbb{C})}=\langle f, (R_0I)^p (g) \rangle_{\mathcal{F}(\mathbb{C})}.$$
We can write $ f(z)=\sum_{n=0}^{\infty} z^n a_n$, $g(z)= \sum_{n=0}^{\infty} z^n b_n$, with $ (a_n)_{n \in \mathbb{N}_0}$, $ (b_n)_{n \in \mathbb{N}_0} \subseteq \mathbb{C}$. From the definition of the operator $ \mathcal{E}_p$ we have that
\begin{equation}
	\label{inn}
\langle \mathcal{E}_p(f), g \rangle_{\mathcal{H}_p(\mathbb{C})}=\left \langle \sum_{n=0}^{\infty} (n+1)^p z^n a_n, g \right \rangle_{\mathcal{F}(\mathbb{C})}=\sum_{n=0}^\infty \frac{n!}{(n+1)^{p}} a_n \overline{b_n}.
\end{equation}
On the other side since
$$ (R_0I)^p z^n= \frac{z^n}{(n+1)^p},$$
we get
\begin{equation}
\label{inn1}
\langle f, (R_0I)^p g \rangle_{\mathcal{F}(\mathbb{C})}=\left \langle f, \sum_{n=0}^{\infty} \frac{z^n a_n}{(n+1)^p} \right \rangle_{\mathcal{H}_p(\mathbb{C})}=\sum_{n=0}^\infty \frac{n!}{(n+1)^{p}} a_n \overline{b_n}.
\end{equation}
Since \eqref{inn} and \eqref{inn1} are equal we get that $ \mathcal{E}_p^*=(R_0I)^p$. Finally, since the operator $\mathcal{E}_p$ is unitary its adjoint coincides with its inverse.
\end{proof}

Now we show that we can write the reproducing kernel of the space $ \mathcal{H}_p(\mathbb{C})$ in terms of the operator $ \mathcal{E}_{p}$.

\begin{lemma}
Let $p \in \mathbb{N}$. Then we have
$$ \mathcal{E}_p^{\bar{w}} \mathcal{E}_p^z (e^{z \bar{w}}):= \frac{e^{z \bar{w}}}{z \bar{w}} T_{2p+1}(z \bar{w}),$$
where $\mathcal{E}_p^{\bar{w}}$, $ \mathcal{E}_p^z$ denote the operator $ \mathcal{E}_p$ with respect to the variables $z$ and $\bar{w}$, respectively.
\end{lemma}
\begin{proof}
Since $\mathcal{E}_p^z (z^n)=(n+1)^p z^n$ and $\mathcal{E}_p^{\bar{w}} (\bar{w}^n)=(n+1)^p \bar{w}^n$ we have
\begin{eqnarray*}
\mathcal{E}_p^{\bar{w}} \mathcal{E}_p^z (e^{z \bar{w}})&=& \mathcal{E}_p^{\bar{w}} \mathcal{E}_p^z \left( \sum_{n=0}^{\infty}  \frac{z^n \bar{w}^n}{n!}\right)\\
&=& \mathcal{E}_p^{\bar{w}} \left( \sum_{n=0}^{\infty}  \frac{(n+1)^p z^n \bar{w}^n}{n!}\right)\\
&=& \sum_{n=0}^{\infty}  \frac{(n+1)^{2p} z^n \bar{w}^n}{n!}\\
&=& K_p(z,w),
\end{eqnarray*}
where $K_p(z,w)$ is the reproducing kernel of the space $ \mathcal{H}_p(\mathbb{C})$. By Proposition \ref{rel} we get the result.
\end{proof}

By combining the properties of the operator $ \mathcal{E}_p$ and the ones of the classical Bargmann transform we can deduce the main properties of the $ \mathcal{H}_p$-Bargmann transform.

\begin{proposition}
Let $p \in \mathbb{N}$. Then the $\mathcal{H}_p$-Bargmann transform is a unitary operator and satisfies
$$ \| \mathcal{B}_p(\varphi) \|_{\mathcal{H}_p(\mathbb{C})}= \| \varphi \|_{L^2(\mathbb{R})}, \qquad \hbox{and} \qquad \mathcal{B}_p(\xi_n)(z)=\frac{(n+1)^p z^n}{\sqrt{n!}}.$$
\end{proposition}
\begin{proof}
First we show that $\mathcal{B}_p$ is an unitary operator. By Proposition \ref{barp}, Theorem \ref{uni} and the fact that the classical Bargmann transform is unitary we get
$$\| \mathcal{B}_p(\varphi) \|_{\mathcal{H}_p(\mathbb{C})}=\| \mathcal{E}_p(\mathcal{B})(\varphi)\|_{\mathcal{H}_p(\mathbb{C})}=\| \mathcal{B}(\varphi) \|_{\mathcal{F}(\mathbb{C})}=\| \varphi \|_{L^2(\mathbb{R})}.$$
By using another time Proposition \ref{barp} and formula \eqref{prop1} we get
$$ \mathcal{B}_p(\xi_n)= \mathcal{E}_p(\mathcal{B}(\xi_n)(z))=\frac{\mathcal{E}_p(z^n)}{\sqrt{n!}}=\frac{(n+1)^p}{\sqrt{n!}} z^n.$$
\end{proof}

\section{The space $\mathcal{F}_p$ and associated operators}

\subsection{Construction and basic properties of the $\mathcal{F}_p$-space}

In this section we study the range of the Schwartz space $ \mathcal{S}(\mathbb{R})$ under the Bargmann transform. The goal of this section is to show the main properties of the space that comes out from this application.

\begin{lemma}
	For a fixed $p  \in \mathbb{N}$ we have
	$$ \mathcal{B} \left(\ell^2(\mathbb{N}, (n+1)^{2p})\right)=\mathcal{F}_p(\mathbb{C}),$$
	where
	\begin{equation}
		\label{space1}
		\mathcal{F}_p(\mathbb{C}):= \left\{ f(z)=\sum_{n=0}^\infty z^n a_n; \quad (a_n)_{n \in \mathbb{N}_0} \subseteq \mathbb{C}, \quad \sum_{n=0}^\infty |a_n|^2 n!(n+1)^{2p}   <\infty \right\}.
	\end{equation}
\end{lemma}
\begin{proof}
The result follows from similar computations performed in Lemma \ref{star}.
\end{proof}

\begin{theorem}
Let $p \in \mathbb{N}$. Then we have
$$ \mathcal{B} \left(\bigcap_{p>0} \ell^2(\mathbb{N}, (n+1)^{2p})\right)= \bigcap_{p>0} \mathcal{B}\left(\ell^2(\mathbb{N}, (n+1)^{2p})\right).$$
\end{theorem}
\begin{proof}
We show the equality by a double inclusion. Let $f \in  \bigcap_{p>0} \ell^2(\mathbb{N}, (n+1)^{2p})$. Then for any $p \in \mathbb{N}$ we have $ f \in \ell^2(\mathbb{N}, (n+1)^{2p})$. Since the Bargmann transform is unitary we have 
$$ \mathcal{B}f \in  \mathcal{B}\left(\ell^2(\mathbb{N}, (n+1)^{2p})\right).$$
This implies that $ \mathcal{B}f \in \bigcap_{p>0} \mathcal{B}\left(\ell^2(\mathbb{N}, (n+1)^{2p})\right).$ 
\\Now we show the other inclusion. Let $ g \in \bigcap_{p>0} \mathcal{B}\left(\ell^2(\mathbb{N}, (n+1)^{2p})\right)$. Then for any $p \in \mathbb{N}$ we have $ g \in  \mathcal{B}\left(\ell^2(\mathbb{N}, (n+1)^{2p})\right)$. Then there exists $\alpha_p \in \ell^2(\mathbb{N}, (n+1)^{2p})$ such that
$$g= \mathcal{B}(\alpha_p).$$
Similarly we can take $ \alpha_q$ such that we have
$$ g= \mathcal{B}(\alpha_q).$$
Therefore we have
$$ \mathcal{B}(\alpha_p)=\mathcal{B}(\alpha_q).$$
Since the Bargmann transform is unitary, for any $p>0$ we have
$$ \alpha_q=\alpha_p:= \alpha \in \ell^2(\mathbb{N}, (n+1)^{2p}).$$
Therefore we have
$$ g = \mathcal{B}(\alpha) \in \mathcal{B} \left(\bigcap_{p>0} \ell^2(\mathbb{N}, (n+1)^{2p})\right).$$
This concludes the proof.
\end{proof}
Now we give a rigorous definition of the space $ \mathcal{F}_p(\mathbb{C})$.
\begin{definition}
	Let $p \in \mathbb{N}$. We define the Hilbert space $ \mathcal{F}_p(\mathbb{C})$ as the set of all entire functions endowed with the inner product given by
	\begin{equation}
		\label{inner01}
		\langle f, g \rangle_{\mathcal{F}_p(\mathbb{C})}= \sum_{n=0}^\infty n!(n+1)^{2p} a_n \overline{b_n},
	\end{equation}
	where $ f(z)= \sum_{n=0}^\infty z^n a_n$ and $g(z)= \sum_{n=0}^\infty z^n b_n $ with $(a_n)_{n \in \mathbb{N}_0}$, $(b_n)_{n \in \mathbb{N}_0} \subseteq \mathbb{C}$.
\end{definition}
The inner product defined in \eqref{inner01} induces the following norm
$$ \| f \|_{\mathcal{F}_p(\mathbb{C})}^2= \sum_{n=0}^\infty n!(n+1)^{2p} | a_n|^2.$$
As an obvious consequence of \eqref{inner01} we have the following
\begin{equation}
	\label{inner1}
	\langle z^n, z^m \rangle_{\mathcal{F}_p(\mathbb{C})}= n! (n+1)^{2p}  \delta_{n,m}.
\end{equation}

\begin{remark}
	For $p=0$ the space $ \mathcal{F}_p(\mathbb{C})$ coincides with the classic Fock space, see \eqref{seq }.
\end{remark}
Let $z$, $w \in \mathbb{C}$. We consider the following function
$$ \mathcal{K}_p(z,w)=\sum_{n=0}^\infty \frac{z^n \bar{w}^n }{n! (n+1)^{2p}}.$$
By using similar arguments of Proposition \ref{necc} we have the following result.
\begin{proposition}
Let $w \in \mathbb{C}$ and $p \in \mathbb{N}$. 
\begin{itemize}
	\item[i)] The function $\mathcal{K}_p(z,w)$ belongs to $ \mathcal{F}_p(\mathbb{C})$.
	\item[ii)] The evaluation map $e_w(f):=f(w)$ is a continuous linear functional on $ \mathcal{F}_p(\mathbb{C})$. Furthermore, we have
	$$|e_w(f)| \leq \| f\|_{\mathcal{F}_p(\mathbb{C})}\| \mathcal{K}_p\|_{\mathcal{H}_p(\mathbb{C})}, \qquad \forall f \in \mathcal{F}_p(\mathbb{C}).$$
\end{itemize}
\end{proposition}
The above result guarantees the existence of a reproducing kernel of the space $\mathcal{F}_p(\mathbb{C})$.
\begin{theorem}
\label{rep2}
Let $p \in \mathbb{N}$. For any $z$, $w \in \mathbb{C}$ the reproducing kernel of the space $\mathcal{F}_p(\mathbb{C})$ is given by $\mathcal{K}_p(z,w)$	
\end{theorem}
\begin{proof}
The result follows from similar arguments used to prove Theorem \ref{kern2}.
\end{proof}
The reproducing kernel of the space $ \mathcal{F}_p(\mathbb{C})$ can be written in terms of peculiar generalized hypergeometric series, that are defined as
\begin{equation}
\label{hyper}
{ }_qF_q(x_1,...,x_q; y_1,...,y_q; z)= \sum_{n=0}^\infty  \frac{(x_1)_n \cdot... \cdot (x_q)_n}{(y_1)_n \cdot... \cdot (y_q)_n} \frac{z^n}{n!}, \qquad z \in \mathbb{C},
\end{equation}
where $(a)_n$ is the Pochhammer symbol defined for $a>0$ as $(a)_n:= \frac{\Gamma(a+n)}{\Gamma(a)}$. 

\begin{proposition}
\label{Poch}
Let $z$, $w \in \mathbb{C}$. Then for $p \in \mathbb{N}$ the reproducing kernel $ \mathcal{K}_p$ of the space $ \mathcal{F}_p(\mathbb{C})$ can be written as
$$ \mathcal{K}_p(z,w)=  { }_{2p}F_{2p}(1,...,1; 2,...,2; z \bar{w}).$$
\end{proposition}
\begin{proof}
By Theorem \ref{rep2} and the expression of the generalized hypergeometric series we have
\begin{eqnarray*}
\mathcal{K}_p(z,w)&=&\sum_{n=0}^\infty \frac{(z \bar{w})^{n}}{(n+1)^{2p} n!}\\
&=& \sum_{n=0}^{\infty} \frac{(n!)^{2p}}{[(n+1)!]^{2p}} \frac{(z \bar{w})^{n}}{n!}\\
&=&\sum_{n=0}^{\infty} \frac{\left(\Gamma(n+1)\right)^{2p}}{\left(\Gamma(n+2)\right)^{2p}} \frac{(z \bar{w})^{n}}{n!}\\
&=&\sum_{n=0}^{\infty} \frac{[(1)_n]^{2p}}{[(2)_n]^{2p}} \frac{(z \bar{w})^n}{n!}\\
&=& { }_{2p}F_{2p}(1,...,1; 2,...,2; z \bar{w}).
\end{eqnarray*}
\end{proof}
 Unlike to what happens for the space $ \mathcal{H}_p(\mathbb{C})$ studied in Section 3 the space $ \mathcal{F}_p(\mathbb{C})$ do not satisfy any geometric description. This is proved in the following result.

\begin{proposition}
Let $p \in \mathbb{N}$. There is not a radial measure $\sigma(x,y)$ such that
\begin{equation}
\label{star2}
\frac{1}{2 \pi} \int_{\mathbb{C}} z^{n} \bar{z}^m d \sigma(x,y)=n! (n+1)^{2p} \delta_{n,m}.
\end{equation}
\end{proposition}
\begin{proof}
To show this result we use Theorem \ref{hamb}. We set $s_n:=n! (n+1)^{2p}$. The matrix 
$$ S_N:= \left((n+m)!(n+m+1)^{2p}\right)_{m,n \geq 0}^N$$
is not positive semidefinite. Indeed, if we pick $n$, $m=0,1$ and $N=1$ we can write the following $2 \times 2$ matrix
$$ \begin{pmatrix}
1 && 4^p\\
4^p && 2 \cdot 9^p
\end{pmatrix}
.$$
This matrix is not positive semidefinite. Hence by Theorem \ref{hamb} we cannot find a radial measure $\sigma(x,y)$ that has a behaviour like in \eqref{star2}.
\end{proof}

\subsection{Some operators on $ \mathcal{F}_p$-space}

The goal of this section is to study the adjoint of the operators $R_0$, $\partial$, $M_z$ and $I$ in the space $ \mathcal{F}_p(\mathbb{C})$. First we show how the adjoint operators act on the monomials.

\begin{theorem}
\label{adjoint}
Let $p \in \mathbb{N}$. Then for $n \in \mathbb{N}$ we have that 
\begin{itemize}
\item the adjoin of the backward-shift operator is given by
\begin{equation}
\label{ad1}
	R_0^* z^n= \frac{(n+1)^{2p-1}}{(n+2)^{2p}} z^{n+1},
\end{equation}
\item the adjoint of the annihilation operator is given by
\begin{equation}
\label{ad2}
	\partial^* z^n= \frac{(n+1)^{2p}}{(n+2)^{2p}} z^{n+1},	
\end{equation}
\item the adjoint of the momentum operator is given by
\begin{equation}
\label{ad3}
M_z^* z^n= \frac{(n+1)^{2p}}{n^{2p-1}}z^{n-1}, \qquad M_z^*(1)=0,
\end{equation}
\item the adjoint of the integration operator is given by
\begin{equation}
\label{ad4}
I^* z^n=\frac{(n+1)^{2p}}{n^{2p}}z^{n-1}, \qquad I^*(1)=0.
\end{equation}
\end{itemize}
\end{theorem}
\begin{proof}
It follows by similar arguments used to prove Proposition \ref{action2}, Proposition \ref{action3}, Proposition \ref{action4} and Proposition \ref{action5} and by using formula \eqref{inner1}.
\end{proof}
In Section 4 we proved that in the space $ \mathcal{H}_p(\mathbb{C})$  the operators $R_0^*$, $\partial^*_z$, $M_z^*$ and $I^*$ can be written as a composition of the respective adjoint operators in the Fock space and a diagonal operator. In the next results we show that there is a similar behaviour for the operators $R_0^*$, $\partial^*_z$, $M_z^*$ and $I^*$ in the space $ \mathcal{F}_p(\mathbb{C})$.

\begin{theorem}
\label{diagope}
Let $p \in \mathbb{N}$. Then we can write the adjoint of the backward-shift operator in the space $ \mathcal{F}_p(\mathbb{C})$ as
\begin{equation}
\label{n1}
R_0^{*}=I(\mathcal{I}-R_0^2IM_z)^{2p}.
\end{equation}
\end{theorem}
\begin{proof}
We show first the result for monomials $z^n$, with $n \in \mathbb{N}_0$. By the binomial theorem we have
\begin{equation}
	\label{binom}
	I \left( \mathcal{I}-R_0^2IM_z\right)^{2p} z^n= I \left( \sum_{k=0}^{2p }\binom{2p}{k} (-1)^k (R_0^2IM_z)^k\right)z^n.
\end{equation}
Now we observe that for $k >0$ we have
\begin{equation}
	\label{aux1}
	(R_0^2IM_z)^k z^n= \frac{1}{(n+2)^k}z^n.
\end{equation}
By plugging \eqref{aux1} into \eqref{binom} we get 
\begin{eqnarray}
	\nonumber
	I \left( \mathcal{I}-R_0^2IM_z\right)^{2p} z^n&=& I \left( \sum_{k=0}^{2p }\binom{2p}{k} \frac{(-1)^k}{(n+2)^k}\right)z^n\\
	\nonumber
	&=& \frac{1}{n+1}  \left( \sum_{k=0}^{2p }\binom{2p}{k} \frac{(-1)^k}{(n+2)^k}\right)z^{n+1}\\
	\nonumber
	&=& \frac{1}{n+1} \left( \sum_{k=0}^{2p-1 }\binom{2p}{k} \frac{(-1)^k}{(n+2)^k}+ \frac{1}{(n+2)^{2p}}\right)z^{n+1}\\
	\label{aux2}
	&=& \frac{1}{n+1}  \left(  \frac{\sum_{k=0}^{2p-1} \binom{2p}{k} (-1)^k  (n+2)^{2p-k}+ 1}{(n+2)^{2p}}\right)z^{n+1}.
\end{eqnarray}

Now, we change the index of the above sum with $\ell=2p-k$ and we get
\begin{eqnarray}
\nonumber
\sum_{k=0}^{2p-1} \binom{2p}{k}  (-1)^k(n+2)^{2p-k}&=& \sum_{\ell=1}^{2p} \binom{2p}{2p - \ell} (n+2)^{\ell} (-1)^{2p - \ell}\\
	\label{sum12}
&=&\sum_{\ell=1}^{2p} \binom{2p}{\ell} (n+2)^{\ell} (-1)^{\ell}.
\end{eqnarray}
By substituting \eqref{sum12} in \eqref{aux2}, by the binomial theorem and by \eqref{ad1} we get
\begin{eqnarray}
\label{starnew}
	I \left( \mathcal{I}-R_0^2IM_z\right)^{2p} z^n&=& \frac{1}{n+1}  \left(  \frac{\sum_{k=1}^{2p} \binom{2p}{k} (n+2)^{k} (-1)^{k}+ 1}{(n+2)^{2p}}\right)z^{n+1}\\
	\nonumber
	&=& \frac{1}{n+1}  \left(  \frac{\sum_{k=0}^{2p} \binom{2p}{k} (-1)^k (n+2)^{k}}{(n+2)^{2p}}\right)z^{n+1}\\
	\nonumber
	&=& \frac{1}{n+1} \frac{(n+1)^{2p}}{(n+2)^{2p}}z^{n+1}\\
	\nonumber
	&=& \frac{(n+1)^{2p-1}}{(n+2)^{2p}}z^{n+1}\\
	\nonumber
	&=& I \left( \mathcal{I}-R_0^2IM_z\right)^{2p} z^n.
\end{eqnarray}	
Now we show the result for a generic function $f$ in $ \mathcal{F}_p(\mathbb{C})$.
Let $f(z)= \sum_{n=0}^\infty a_n z^n$ and $g(z)= \sum_{n=0}^\infty b_n z^n$ be another function in $ \mathcal{F}_p(\mathbb{C})$, with $ (a_n)_{n \in \mathbb{N}_0}$, $ (b_n)_{n \in \mathbb{N}_0} \subseteq \mathbb{C}$. By the binomial theorem we get
\begingroup\allowdisplaybreaks	
\begin{eqnarray*}
	\langle f,  I \left( \mathcal{I}-R_0^2IM_z\right)^{2p}   g\rangle_{\mathcal{F}_p(\mathbb{C})}&=&\left \langle f,  I \left( \mathcal{I}-R_0^2IM_z\right)^{2p} \sum_{n=0}^\infty b_n z^n\right \rangle_{\mathcal{F}_p(\mathbb{C})}\\
	&=& \left \langle f,   \sum_{n=0}^\infty  \frac{(n+1)^{2p-1}}{(n+2)^{2p}}z^{n+1} b_{n} \right \rangle_{\mathcal{F}_p(\mathbb{C})}\\
	&=& \left \langle f,   \sum_{n=1}^\infty  \frac{n^{2p}}{n(n+1)^{2p}}z^{n}b_{n-1} \right \rangle_{\mathcal{F}_p(\mathbb{C})}\\
	&=&  \sum_{n=1}^\infty  a_{n} \overline{b_{n-1}} (n-1)! n^{2p}\\
	&=&  \left \langle \sum_{n=0}^\infty a_{n+1}z^n, g \right \rangle_{\mathcal{F}_p(\mathbb{C})} \\
	&=& \langle R_0 f, g \rangle_{\mathcal{F}_p(\mathbb{C})}.
\end{eqnarray*}
\endgroup
This shows that $R_0^* = I \left( \mathcal{I}-R_0^2IM_z\right)^{2p} .$
\end{proof}

\begin{theorem}
Let $p \in \mathbb{N}$. Then we can write the adjoint annihilation operator for the space $ \mathcal{F}_p(\mathbb{C})$ as
\begin{equation}
	\label{n2}
	\partial^*= M_z(\mathcal{I}-R_0^2IM_z)^{2p},
\end{equation}
\end{theorem}
\begin{proof}
We show first the result for monomials $z^n$, with $n \in \mathbb{N}_0$. By similar arguments used to prove \eqref{starnew} we deduce that
$$ (\mathcal{I}-R_0^2IM_z)^{2p} z^n= \frac{(n+1)^{2p}}{(n+2)^{2p}}z^{n}.$$
By \eqref{ad2} we get
\begin{equation}
	\label{a21}
M_z(\mathcal{I}-R_0^2IM_z)^{2p} z^n=\frac{(n+1)^{2p}}{(n+2)^{2p}}z^{n+1}=\partial^{*}z^n.
\end{equation}
Now we show the result for a generic function $f$ in $ \mathcal{F}_p(\mathbb{C})$.
Let us consider $f(z)= \sum_{n=0}^\infty a_n z^n$ and $g(z)= \sum_{n=0}^\infty b_n z^n$ in $ \mathcal{F}_p(\mathbb{C})$, with $ (a_n)_{n \in \mathbb{N}_0}$, $ (b_n)_{n \in \mathbb{N}_0} \subseteq \mathbb{C}$. By the binomial theorem we get
\begingroup\allowdisplaybreaks	
\begin{eqnarray*}
	\langle f,  M_z(\mathcal{I}-R_0^2IM_z)^{2p} g\rangle_{\mathcal{F}_p(\mathbb{C})}&=&\left \langle f,  M_z(\mathcal{I}-R_0^2IM_z)^{2p} \sum_{n=0}^\infty b_n z^n\right \rangle_{\mathcal{F}_p(\mathbb{C})}\\
	&=& \left \langle f,   \sum_{n=0}^\infty  \left(\frac{n+1}{n+2}\right)^{2p} z^{n+1}b_n \right \rangle_{\mathcal{F}_p(\mathbb{C})}\\
	&=& \left \langle f,   \sum_{n=1}^\infty  \frac{n^{2p}}{(n+1)^{2p}} z^{n}b_{n-1} \right \rangle_{\mathcal{F}_p(\mathbb{C})}\\
	&=&  \sum_{n=1}^\infty  a_{n} \overline{b_{n-1}} n! n^{2p} \\
	&=&  \left \langle \sum_{n=0}^\infty (n+1)a_{n+1}z^n, g \right \rangle_{\mathcal{F}_p(\mathbb{C})} \\
	&=& \langle \partial f, g \rangle_{\mathcal{F}_p(\mathbb{C})}.
\end{eqnarray*}
\endgroup
This shows that $\partial^* = M_z(\mathcal{I}-R_0^2IM_z)^{2p}$.

\end{proof}

\begin{theorem}
\label{aux12}
Let $p \in \mathbb{N}$. Then we can write  adjoint of the momentum operator for the space $ \mathcal{F}_p(\mathbb{C})$ as	
\begin{equation}
	\label{n3}
	M_z^*= \partial(\mathcal{I}+IR_0)^{2p}
\end{equation}
\end{theorem}
\begin{proof}
We show first the result for monomials $z^n$, with $n \in \mathbb{N}$. By the binomial theorem formula we have
\begin{equation}
	\label{newt}
\partial(\mathcal{I}+IR_0)^{2p} z^n= \partial \left(\sum_{k=0}^{2p} \binom{2p}{k} (IR_0)^{k}\right)z^{n}.
\end{equation}
We observe that for $k \geq 0$ we have
\begin{equation}
	\label{negfund}
	(IR_0)^k z^n= \frac{z^n}{n^k}.
\end{equation}
By inserting \eqref{negfund} into \eqref{newt} we get
\begin{eqnarray}
	\nonumber
	\partial(\mathcal{I}+IR_0)^{2p} z^n&=& \partial \left(\sum_{k=0}^{2p} \binom{2p}{k}   \frac{1}{n^k}\right)z^{n}\\
	\nonumber
	\nonumber
	&=& n \left(\sum_{k=0}^{2p-1} \binom{2p}{k} \frac{1}{n^{k}}+ \frac{1}{n^{2p}} \right) z^{n}\\
	\label{n6}
	&=& n \left(\frac{\sum_{k=0}^{2p-1} \binom{2p}{k} n^{2p-k}+1}{n^{2p}}\right) z^{n}.
\end{eqnarray}
Now, we change the index of the sum above  with $\ell=2p-k$ and we get
\begin{equation}
	\label{change}
	\sum_{k=0}^{2p-1} \binom{2p}{k}  n^{2p-k}= \sum_{\ell=1}^{2p} \binom{2p}{2p - \ell} n^{\ell}.
\end{equation}
We observe that $ \binom{2p}{2p - \ell}= \binom{2p}{\ell}$. This implies that
\begin{equation}
	\label{sum}
\sum_{k=0}^{2p-1} \binom{2p}{k} \frac{n^{2p}}{n^k}+1=\sum_{\ell=1}^{2p} \binom{2p}{2p - \ell} n^{\ell}+1=\sum_{\ell=0}^{2p} \binom{2p}{\ell} n^{\ell}.
\end{equation}
Finally by inserting \eqref{sum} in \eqref{n6} we get
\begin{equation}
\label{aux}
\partial(\mathcal{I}+IR_0)^{2p} z^n=\partial \left(\frac{\sum_{\ell=0}^{2p} \binom{2p}{\ell} n^{\ell}}{n^{2p}}\right) z^n=\frac{(n+1)^{2p}}{n^{2p-1}}z^{n-1}.
\end{equation}
Finally by \eqref{ad3} we have
\begin{equation}
\label{ns}
\partial(\mathcal{I}+IR_0)^{2p} z^n= M_{z}^{*}z^n.
\end{equation}
Formula \eqref{ns} is still valid for $n=0$. Since 
$$ M_z^*(1)=0, \qquad \partial(\mathcal{I}+IR_0)^{2p}(1)=0.$$
Now, we prove the result for a generic $f$ in $\mathcal{F}_p(\mathbb{C})$.
Let $f(z)= \sum_{n=0}^\infty a_n z^n$ and $g(z)= \sum_{n=0}^\infty b_n z^n$ be another function in $ \mathcal{F}_p(\mathbb{C})$, with $ (a_n)_{n \in \mathbb{N}_0}$, $ (b_n)_{n \in \mathbb{N}_0} \subseteq \mathbb{C}$, then we have
\begingroup\allowdisplaybreaks	
\begin{eqnarray*}
	\langle f,  \partial(\mathcal{I}+IR_0)^{2p} g\rangle_{\mathcal{F}_p(\mathbb{C})}&=&\left \langle f,  \partial(\mathcal{I}+IR_0)^{2p} \sum_{n=0}^\infty b_n z^n\right \rangle_{\mathcal{F}_p(\mathbb{C})}\\
		&=& \left \langle f,   \sum_{n=1}^\infty  \frac{(n+1)^{2p}}{n^{2p-1}} z^{n-1}b_n \right \rangle_{\mathcal{F}_p(\mathbb{C})}\\
	&=& \left \langle f,   \sum_{n=0}^\infty  \frac{(n+2)^{2p}}{(n+1)^{2p-1}} z^{n}b_{n+1} \right \rangle_{\mathcal{F}_p(\mathbb{C})}\\
	&=&  \sum_{n=0}^\infty  a_{n} \overline{b_{n+1}}(n+1)!(n+2)^{2p} \\
	&=&  \left \langle \sum_{n=1}^\infty a_{n-1}z^n, g \right \rangle_{\mathcal{F}_p(\mathbb{C})} \\
	&=& \langle M_z f, g \rangle_{\mathcal{F}_p(\mathbb{C})}.
\end{eqnarray*}
\endgroup
This shows that $M_z^* = \partial(\mathcal{I}+IR_0)^{2p}$.

\end{proof}

\begin{theorem}
Let $p \in \mathbb{N}$. Then we can write the adjoint  of the integration operator for the space $ \mathcal{F}_p(\mathbb{C})$ as
\begin{equation}
\label{n4}
I^*= R_0(\mathcal{I}+IR_0)^{2p}
\end{equation}
\end{theorem}
\begin{proof}
We show first the result for monomials $z^n$ with $n \in \mathbb{N}$. By similar arguments used to prove \eqref{aux} we have

$$ (\mathcal{I}+IR_0)^{2p} z^n=  \frac{(n+1)^{2p}}{n^{2p}} z^n.$$
Hence by \eqref{ad4} we have
$$ R_0 (\mathcal{I}+IR_0)^{2p} z^n=\frac{(n+1)^{2p}}{n^{2p}} z^{n-1}=I^* z^n.$$
If we consider $n=0$ we have
$$ I^*(1)=0, \qquad R_0 (\mathcal{I}+IR_0)^{2p} (1)=0.$$
Now, we prove the result for a generic $f \in \mathcal{F}_p(\mathbb{C})$.
Let $f(z)= \sum_{n=0}^\infty a_n z^n$ and $g(z)= \sum_{n=0}^\infty b_n z^n$ in $ \mathcal{F}_p(\mathbb{C})$, with $ (a_n)_{n \in \mathbb{N}_0}$, $ (b_n)_{n \in \mathbb{N}_0} \subseteq \mathbb{C}$. By the binomial theorem we get
\begingroup\allowdisplaybreaks	
\begin{eqnarray*}
	\langle f,  R_0 (\mathcal{I}+IR_0)^{2p}g\rangle_{\mathcal{F}_p(\mathbb{C})}&=&\left \langle f,  \partial (\mathcal{I}+IR_0)^{2p} \sum_{n=0}^\infty b_n z^n\right \rangle_{\mathcal{F}_p(\mathbb{C})}\\
	&=& \left \langle f,   \sum_{n=1}^\infty  \frac{(n+1)^{2p}}{n^{2p}} z^{n-1}b_n \right \rangle_{\mathcal{F}_p(\mathbb{C})}\\
	&=& \left \langle f,   \sum_{n=0}^\infty  \frac{(n+2)^{2p}}{(n+1)^{2p}} z^{n}b_{n+1} \right \rangle_{\mathcal{F}_p(\mathbb{C})}\\
	&=&  \sum_{n=0}^\infty  a_{n} \overline{b_{n+1}} n!(n+2)^{2p} \\
	&=&  \left \langle \sum_{n=1}^\infty \frac{a_{n-1}}{n}z^n, g \right \rangle_{\mathcal{F}_p(\mathbb{C})} \\
	&=& \langle I f, g \rangle_{\mathcal{F}_p(\mathbb{C})}.
\end{eqnarray*}
\endgroup
This shows that $I^{*} = R_{0}(\mathcal{I}+IR_0)^{2p}$.
 
\end{proof}

We sum up all the previous results in the following table

\begin{center}
	\begin{tabular}{| l | l |l |l|l|}
		\hline
		\rule[-4mm]{0mm}{1cm}
		{\bf Operators} & {$R_0^{*}$} & {$\partial^*$ } & {$M_z^{*}$} &{$I^{*}$}\\
		\hline
		\rule[-4mm]{0mm}{1cm}
		{} & $I(\mathcal{I}-R_0^2IM_z)^{2p}$& $M_z(\mathcal{I}-R_0^2IM_z)^{2p}$ & $ \partial  (\mathcal{I}+IR_0)^{2p}$& $ R_{0}(\mathcal{I}+IR_0)^{2p}$ \\
		\hline
		\hline
	\end{tabular}
\end{center}

Now we investigate the commutators of $M_z$ and $M_z^*$ as well as $R_0$ and $R_0^{*}$. Similarly to what happens for the space $ \mathcal{H}_p(\mathbb{C})$ we do not have any commutation relations, like in the classic case. In the next result we use the notation of backwards of a diagonal operator, see \eqref{back}.
\begin{proposition}
Let $p \in \mathbb{N}$. The operators $M_z$ and $M_z^*$ satisfy the following formal relation
\begin{equation}
\label{ncomm}
[M_z, M_z^*]=\mathsf{A}(Id+\mathsf{E})^{2p}- \mathsf{A}^{(-1)}(Id+[\mathsf{E}^{(-1)}])^{2p},
\end{equation}
where $Id$ is the identity matrix, $ \mathsf{A}:= \hbox{diag}(0,1,2,3,...)$ and $\mathsf{E}:= \hbox{diag}\left(1,1, \frac{1}{2}, \frac{1}{3},..., \frac{1}{n},...\right)$.
\end{proposition}
\begin{proof}
We show the theorem only for monomials $z^n$, with $n \in \mathbb{N}$. By Theorem \ref{aux12} we have
\begin{eqnarray*}
[M_z, M_z^*](z^n)&=& M_z M_z^*(z^n)-M_z^*M_z(z^n)\\
&=& \sum_{k=0}^{2p} \binom{2p}{k}M_z \partial (IR_0)^k (z^n)- \sum_{k=0}^{2p} \binom{2p}{k} \partial (IR_0)^{k}(z^{n+1})\\
&=& \left[n \sum_{k=0}^{2p} \binom{2p}{k} \frac{1}{n^k}-(n+1) \sum_{k=0}^{2p}\binom{2p}{k} \frac{1}{(n+1)^k}\right]z^n.
\end{eqnarray*}
By the binomial formula we have
$$ [M_z^, M_z^*](z^n)= \left[n \left(1+ \frac{1}{n}\right)^{2p}-(n+1) \left(1+ \frac{1}{n+1}\right)^{2p}\right]z^n.$$
If we take $n=0$ we have
$$ [M_z, M_z^*](1)=M_z M_z^*(1)-M_z^*M_z(1)=-4^p.$$
Therefore
$$ [M_z^, M_z^*](z^n)=\begin{cases}
-4^p, \qquad n=0\\	
\left[n \left(1+ \frac{1}{n}\right)^{2p}-(n+1) \left(1+ \frac{1}{n+1}\right)^{2p}\right]z^n, \quad n=1,2,3,..
\end{cases}
$$
Finally by the notion of the backwards of a diagonal operator, see  \eqref{back}, we have
$$[M_z, M_z^*]= \left[\mathsf{A}(Id+\mathsf{E})^{2p}- \mathsf{A}^{(-1)}(Id+[\mathsf{E}^{(-1)}])^{2p}\right].$$
\end{proof}

\begin{proposition}
Let $p \in \mathbb{N}$. The operators $R_0$ and $R_0^*$ satisfy the following formal relation
\begin{equation}
\label{ncomm1}
[R_0, R_0^*]=D_0 \left[\mathsf{A}^{(-1)} \left(Id-[\mathsf{A}^{(-1)}]^{-1}\right)^{2p}-\mathsf{A} \left(Id-[\mathsf{A}^{(-2)}]^{-1}\right)^{2p}  \right].,
\end{equation}
where $Id$ is the identity matrix, $ \mathsf{A}:= \hbox{diag}(0,1,2,3,...)$, $[.]^{-1}$ is the inverse matrix and $D_0$ is a diagonal matrix introduced in in \eqref{mat}.
\end{proposition}
\begin{proof}
We prove the result for monomials $z^n$, for $n \in \mathbb{N}$. By Theorem \ref{diagope} we have
\begin{eqnarray*}
[R_0, R_0^*](z^n)&=& R_0 R_0^*(z^n)-R_0^{*}R_0(z^n)\\
&=& \sum_{k=0}^{2p} \binom{2p}{k} (-1)^k R_0I (R_0^2IM_z)(z^n)- \sum_{k=0}^{2p} \binom{2p}{k}(-1)^k I(R_0^2IM_z)(z^{n-1})\\
&=& \left(\sum_{k=0}^{2p} \binom{2p}{k} (-1)^k \frac{1}{(n+2)^k (n+1)}- \sum_{k=0}^{2p} \binom{2p}{k} (-1)^k \frac{1}{(n+1)^k n}\right)z^n.
\end{eqnarray*}
By the binomial formula we have
\begin{eqnarray*}
[R_0, R_0^*](z^n)&=& \left[\frac{1}{n+1} \left(1- \frac{1}{n+2}\right)^{2p}- \frac{1}{n} \left(1- \frac{1}{n+1}\right)^{2p}\right]z^n\\
&=&- \frac{1}{n(n+1)}\left[(n+1) \left(1- \frac{1}{n+1}\right)^{2p}- n \left(1- \frac{1}{n+2}\right)^{2p}\right]z^n.
\end{eqnarray*}

For $n=0$ we have
$$[R_0, R_0^*](1)=R_0 R_0^*(1)-R_0^{*}R_0(1)=\frac{1}{4^p}.$$
Therefore
$$ [R_0, R_0^*](z^n)=\begin{cases}
\frac{1}{4^p}, \qquad n=0\\
- \frac{1}{n(n+1)}\left[(n+1) \left(1- \frac{1}{n+1}\right)^{2p}- n \left(1- \frac{1}{n+2}\right)^{2p}\right]z^n, \qquad n=1,2,...
\end{cases}
$$
Finally by the notion of the backwards of a diagonal operator, see  \eqref{back}, and the notation \eqref{mat} we have
$$[R_0, R_0^*]=D_0 \left[\mathsf{A}^{(-1)} \left(Id-[\mathsf{A}^{(-1)}]^{-1}\right)^{2p}-\mathsf{A} \left(Id-[\mathsf{A}^{(-2)}]^{-1}\right)^{2p}  \right].$$
\end{proof}

\begin{remark}
By taking $p=0$ in \eqref{ncomm} and \eqref{ncomm1} we get back the relations \eqref{comm} and \eqref{com}, respectively. Indeed
$$D_0 \left[\mathsf{A}^{(-1)}- \mathsf{A}\right]=D_0 \left[\hbox{diag}\left(1, 2, 3,...\right)-\hbox{diag}\left(0,1, 2,3,...\right)\right]=D_0Id =D_0.$$
\end{remark}

\subsection{The $\mathcal{F}_p$-Bargmann transform}

The goal of this section is to develop a further characterization of the space $ \mathcal{F}_p(\mathbb{C})$ through a unitary operator from $L^2(\mathbb{R})$ onto $ \mathcal{F}_p(\mathbb{C})$. We find out that this operator can be considered a generalization of the Bargmann transform in the space $ \mathcal{F}_p(\mathbb{C})$.
\begin{proposition}
	Let $p \in \mathbb{N}$. For every $x \in \mathbb{R}$ and $z \in \mathbb{C}$ we define the function
	\begin{equation}
		\mathcal{A}_{p}(z,x)=\sum_{n=0}^{\infty} \frac{z^n}{(n+1)^{p}  \sqrt{n!}} \xi_n(x),
	\end{equation}
	where $\xi_n$ are the normalized Hermite  functions. Then we have
	\begin{itemize}
		\item[1)] The function $	\mathcal{A}_p(.,x)$ is entire for every $x \in \mathbb{R}$.
		\item[2)] A function $f \in \mathcal{F}_p(\mathbb{C})$ if and only if there exists $\varphi \in L^{2}(\mathbb{R})$ such that
		\begin{equation}
			f(z)= \int_{\mathbb{R}} 	\mathcal{A}_p(z,x) \varphi(x) dx= \langle \varphi, \overline{	\mathcal{A}_p(z,.)} \rangle_{L^2}.
		\end{equation}
	\end{itemize}
\end{proposition}
\begin{proof}
The result follows from similar arguments used to prove Proposition \ref{kernel1}.
\end{proof}

The above result motivates the following
\begin{definition}
\label{bargp}
For $p \in \mathbb{N}$ and any $\varphi \in L^2(\mathbb{R})$ we define the $\mathcal{F}_p$-Bargmann transform as
$$ \mathcal{SB}_p(\varphi(z)):= \int_{\mathbb{R}} \mathcal{A}_p(z,x) \varphi(x) dx, \qquad \mathcal{A}_p(z,x):= \frac{1}{(2 \pi)^{\frac{1}{4}}} \sum_{n=0}^{\infty} \frac{z^n h_n(x)}{(n+1)^p n! }.$$
\end{definition}

By using the notations fixed in \eqref{nota1} we have

\begin{proposition}
Let $p \in \mathbb{N}$. For any $z$, $w \in \mathbb{C}$ we have
$$ \langle \mathcal{A}_p^z, \mathcal{A}_p^w \rangle_{L^2}={ }_{2p}F_{2p}(1,...,1; 2,...,2; z \bar{w}),$$
where ${ }_{2p}F_{2p}$ is the generalized hypergeometric series.
\end{proposition}
\begin{proof}
By the definition of the kernel $ \mathcal{A}_p$ and Proposition \ref{Poch} we have
\begin{eqnarray*}
\langle \mathcal{A}_p^z, \mathcal{A}_p^w \rangle_{L^2}&=& \int_{\mathbb{R}} \left( \sum_{j=0}^{\infty}\frac{z^j \xi_j(x)}{(j+1)^p\sqrt{j!}}\right) \left( \sum_{k=0}^{\infty}\frac{\bar{w}^k \xi_k(x)}{(k+1)^p\sqrt{k!}}\right)dx\\
&=& \sum_{j,k=0}^\infty  \frac{z^j \bar{w}^k}{(j+1)^p (k+1)^p \sqrt{j!} \sqrt{k!}} \int_{\mathbb{R}} \xi_j(x) \xi_k(x) dx\\
&=& \sum_{j=0}^{\infty} \frac{z^j \bar{w}^j}{j! (j+1)^{2p}}\\
&=& \mathcal{K}_p(z,w)\\
&=& { }_{2p}F_{2p}(1,...,1; 2,...,2; z \bar{w}).
\end{eqnarray*}
\end{proof}

The kernel of the $\mathcal{F}_p$-Bargmann transform can be written in a peculiar way. 

\begin{theorem}
\label{ker4}
For $p \in \mathbb{N}$ we have
\begin{equation}
\mathcal{A}_p(z,x)= (R_0I)^p A(z,x), \qquad \forall (z,x) \in \mathbb{C} \times \mathbb{R},
\end{equation}
where $A(z,x)$ is the kernel of the $\mathcal{F}_p$-Bargmann transform.
\end{theorem}
\begin{proof}
We observe that $(R_0I)^p= \frac{z^n}{(n+1)^p}$ then by the definition of the kernel of the classical Bargmann we have
\begin{eqnarray*}
\mathcal{A}_p(z,x)&=& \frac{1}{(2 \pi)^{\frac{1}{4}}}\sum_{n=0}^{\infty} \frac{z^n h_n(x)}{(n+1)^p n! 2^{\frac{n}{2}}}\\
&=& (R_0I)^p \left(\frac{1}{(2 \pi)^{\frac{1}{4}}}\sum_{n=0}^{\infty} \frac{z^n h_n(x)}{n! 2^{\frac{n}{2}}}\right)\\
&=&(R_0I)^p A(z,x).
\end{eqnarray*}
\end{proof}
The following result, see \cite[Thm. 4.6]{ACK}, will be useful to write the kernel $ \mathcal{A}_p(z,x)$ in another way.
\begin{proposition}
\label{rel1}
Given $D_0$ as \eqref{mat} we have
\begin{equation}
\label{app3}
(IR_0)^n=\sum_{k=1}^{n} \Lambda_{k,n} I^k R_0^k,
\end{equation}
where $ \Lambda_{k,n}$ is a diagonal operator defined by means of the following recurrence relation
\begin{equation}
\label{coeff}
\Lambda_{k,n}:= \Lambda_{k-1,n-1}+\sum_{\ell=1}^{k} D_0^{(\ell)}\Lambda_{k,n-1},
\end{equation}
where $D_0^{(\ell)}$ is the forward diagonal of $D_0$, see \eqref{forw}. 
The initial values are given by
$$ \Lambda_{n,n}= Id \qquad \hbox{and} \qquad \Lambda_{0,n}=\hbox{diag}(0,0,0,...), \quad n \in \mathbb{N},$$
and for the remaining values $k$, $p$ we have
$$\Lambda_{k,n}=\hbox{diag}(0,0,0,...).$$
\end{proposition}

\begin{proposition}
\label{ker3}
Let $p \in \mathbb{N}$ then for any $z$, $w \in \mathbb{C}$ we have
\begin{eqnarray}
\label{ker2}
\mathcal{A}_p(z,w)&=& (D_0+IR_0) ^pA(z,w)\\
\label{ker1}
&=& \sum_{j=0}^{p} \sum_{k=1}^{j}D_{0}^{p-j} \binom{p}{j} \Lambda_{k,j} I^k R_0^k A(z,w),
\end{eqnarray}
where $\Lambda_{k,j}$ are defined in \eqref{coeff}.
\end{proposition}
\begin{proof}
The equality \eqref{ker2} follows by Theorem \ref{ker4} and the relation \eqref{com}. Now we show \eqref{ker1}. By the binomial formula and Proposition \ref{rel1} we have
\begin{eqnarray*}
\mathcal{A}_p(z,w)&=& (R_0I)^p A(z,w)\\
&=& \sum_{j=0}^{p} \binom{p}{j} D_{0}^{p-j} (IR_0)^j A(z,w)\\
&=&\sum_{j=0}^{p} \sum_{k=1}^{j}D_{0}^{p-j} \binom{p}{j} \Lambda_{k,j} I^k R_0^k A(z,w).
\end{eqnarray*}
\end{proof}

The previous result paved the way to write the $\mathcal{F}_p$-Bargmann transform in three different ways.

\begin{proposition}
\label{opBarg}
For $p \in \mathbb{N}$ and any $ \varphi \in L^2(\mathbb{R})$ we can write the $\mathcal{F}_p$-Bargmann transform as
\begin{eqnarray*}
\mathcal{SB}_p(\varphi(z))&=&(R_0I)^p (\mathcal{B} \varphi(z))\\
&=& (D_0+IR_0)^p(\mathcal{B} \varphi(z))\\
&=&\sum_{j=0}^{p} \sum_{k=1}^{j}D_{0}^{p-j} \binom{p}{j} \Lambda_{k,j} I^k R_0^k (\mathcal{B} \varphi(z)),
\end{eqnarray*}
where $D_0$ is the diagonal matrix defined in \eqref{mat}.
\end{proposition}
\begin{proof}
The result follows by Definition \ref{bargp}, Theorem \ref{ker4} and Proposition \ref{ker3}.
\end{proof}

\subsection{The Fock space connected to the $ \mathcal{F}_p$-space}

In this Section we show a connection between the Fock space and $ \mathcal{F}_p(\mathbb{C})$. Similarly to what happens for the space $ \mathcal{H}_p(\mathbb{C})$ in Section 6 the connection is provided by a diagonal operator. Moreover, in this section we further investigate the properties of the $\mathcal{F}_p$-Bargmann transform.
\begin{equation}
\label{ope5}
\mathcal{V}_p:=(R_0I)^p.
\end{equation}

\begin{theorem}
\label{uni2}
Let $p \in \mathbb{N}$. Then the operator $\mathcal{V}_p$ is a surjective isometry from the space $ \mathcal{F}(\mathbb{C})$ onto $ \mathcal{F}_p(\mathbb{C})$. Moreover the adjoint of $\mathcal{V}_p$ is given by the operator $ \mathcal{E}_p$  defined in \eqref{ope2}.
\end{theorem}
\begin{proof}
Firstly we show that the operator $ \mathcal{V}_p$ is an isometric operator. We observe that $ \mathcal{V}_p(z^n)= \frac{z^n}{(n+1)^p}$. Therefore for $f \in \mathcal{F}(\mathbb{C})$, such that $f(z)= \sum_{n=0}^{\infty} z^n a_n$, with $(a_n)_{n \in \mathbb{N}_0} \subseteq \mathbb{C}$, we have
$$
\mathcal{V}_p(f)= \sum_{n=0}^\infty \frac{z^n}{(n+1)^p} a_n= \sum_{n=0}^\infty z^n \beta_n,
$$
where $ \beta_n:= \frac{a_n}{(n+1)^p}$. Hence we get
$$ \| \mathcal{V}_p(f) \|_{\mathcal{F}_p(\mathbb{C})}^2= \sum_{n=0}^{\infty} n! (n+1)^{2p} |\beta_n|^2= \|f \|_{\mathcal{F}(\mathbb{C})}^2.$$
Now, we show that the operator $ \mathcal{V}_p$ is surjective. Let $g \in \mathcal{F}_p(\mathbb{C})$. Hence we have $g(z)= \sum_{n=0}^{\infty} z^n a_n$ with $(a_n)_{n \in \mathbb{N}_0} \subseteq \mathbb{C}$. Our goal is to find a  function $g \in \mathcal{F}(\mathbb{C})$ such that
\begin{equation}
\label{rela}
g(z)= \mathcal{V}_p(f)(z).
\end{equation}
We have to choose a sequence $ (b_n)_{n \in \mathbb{N}_0} \subseteq \mathbb{C}$ such that $g(z)= \sum_{n=0}^{\infty} z^n b_n$ and $\| g \|^2_{\mathcal{F}(\mathbb{C})}< \infty$. If we take $b_n:= (n+1)^p a_n$, we have
$$
\mathcal{V}_p(f)= \sum_{n=0}^{\infty} \frac{z^n}{(n+1)^p}b_n= g(z).
$$
Now, since $f \in \mathcal{F}_p(\mathbb{C})$ we have
$$ \| g\|_{\mathcal{F}(\mathbb{C})}^2=\sum_{n=0}^{\infty} n! |b_n|^2=\sum_{n=0}^{\infty} n! (n+1)^{2p}|a_n|^2< \infty.$$
Therefore $g \in \mathcal{F}(\mathbb{C})$. 
\\Now we show that the adjoint of $ \mathcal{V}_p$ is $ \mathcal{E}_p$. We have to show that
$$ \langle \mathcal{V}_p(f),g \rangle_{\mathcal{F}_p(\mathbb{C}) }= \langle f, \mathcal{E}_p(g) \rangle_{\mathcal{F}(\mathbb{C})},$$
for $f \in \mathcal{F}(\mathbb{C})$ and $g \in \mathcal{F}_p(\mathbb{C})$. Let $ f(z)=\sum_{n=0}^{\infty} z^n a_n$, $g(z)= \sum_{n=0}^{\infty} z^n b_n$, with $ (a_n)_{n \in \mathbb{N}_0}$, $ (b_n)_{n \in \mathbb{N}_0} \subseteq \mathbb{C}$ then we have
\begin{equation}
\label{adj}
\langle \mathcal{V}_p(f),g \rangle_{\mathcal{F}_p(\mathbb{C}) }=\sum_{n=0}^{\infty} n! (n+1)^p a_n \overline{b_n}.
\end{equation}
We observe that $ \mathcal{E}_p(z^n)=(n+1)^p z^n$. This implies that
\begin{equation}
\label{adj2}
\langle f, \mathcal{E}_p(g) \rangle_{\mathcal{F}(\mathbb{C})}=\sum_{n=0}^{\infty} n! (n+1)^p a_n \overline{b_n}.
\end{equation}
Since \eqref{adj} and \eqref{adj2} are equal we get the result.
\end{proof}

\begin{lemma}
Let $p \in\mathbb{N}$. Then for any $z$ $w \in\mathbb{C}$ we have
$$ \mathcal{V}_{p}^{\bar{w}}\mathcal{V}_{p}^{z}(e^{z \bar{w}})={ }_{2p}F_{2p}(1,...,1; 2,...,2; z \bar{w}),$$
where $\mathcal{V}_{p}^{\bar{w}}$ $\mathcal{V}_{p}^{z}$ denote the application of the operator $ \mathcal{V}_p$ with respect to the variables $\bar{w}$ and $z$, respectively.
\end{lemma}
\begin{proof}
We start by observing that $ \mathcal{V}_{p}^{\bar{w}}(\bar{w}^n)= \frac{\bar{w}^n}{(n+1)^p}$ and $ \mathcal{V}_{p}^{z}(z^n)= \frac{z^n}{(n+1)^p}$. By Proposition \ref{Poch} we have
	\begin{eqnarray*}
	\mathcal{V}_{p}^{\bar{w}}\mathcal{V}_{p}^{z}(e^{z \bar{w}})&=& \sum_{n=0}^{\infty} \frac{z^n \bar{w}^n}{n! (n+1)^{2p}}\\
	&=& \mathcal{K}_p(z,w)\\
	&=& { }_{2p}F_{2p}(1,...,1; 2,...,2; z \bar{w}).
\end{eqnarray*}
\end{proof}

\begin{proposition}
Let $p \in \mathbb{N}$. Then the $\mathcal{F}_p$-Bargmann transform is a unitary operator and holds that
$$ \| \mathcal{SB}_p(\varphi)\|_{\mathcal{F}_p(\mathbb{C})}=\| \varphi \|_{L^2(\mathbb{R})}, \qquad \hbox{and} \qquad \mathcal{B}_p(\xi_n)(z)= \frac{z^n}{(n+1)^p \sqrt{n!}}.$$
\end{proposition}
\begin{proof}
Since the Bargmann transform  $\mathcal{SB}_p$ is a composition of two unitary operators: the classical Bargmann transform and the operator $\mathcal{V}_p$, see Theorem \ref{uni2}, we get 
$$ \| \mathcal{SB}_p(\varphi)\|_{\mathcal{F}_p(\mathbb{C})}=\| \varphi \|_{L^2(\mathbb{R})}.$$
Finally by formula \eqref{prop1}, the fact that $ \mathcal{V}_p(z^n)= \frac{z^n}{(n+1)^p}$ and by using another time Proposition \ref{opBarg} we get that
$$ \mathcal{SB}_p(\xi_n)(z)= \frac{z^n}{(n+1)^p \sqrt{n!}}.$$
\end{proof}

\section{A relation between the $\mathcal{H}_p$ and $\mathcal{F}_p$ spaces}

In this section we show how the two spaces $ \mathcal{H}_p(\mathbb{C})$ and $ \mathcal{F}_p(\mathbb{C})$ are connected each other through the Fock space.

\begin{theorem}
\label{conn}
Let $p \in \mathbb{N}$. The operator that maps the space $ \mathcal{F}_p(\mathbb{C}) $ onto $ \mathcal{H}_p(\mathbb{C})$ is given by
$$ \Theta_p= \left(\mathcal{I}+M_z \partial\right)^{2p}.$$
\end{theorem}
\begin{proof}
By combining Theorem \ref{uni} and Theorem \ref{uni2} we have the following scheme

\[
\begin{tikzcd}
	& \mathcal{F}_{p}(\mathbb{C})  \arrow{dr}{\Theta_p} \\
	\mathcal{F}(\mathbb{C}) \arrow{ur}{\mathcal{V}_p} \arrow{rr}{\mathcal{E}_p} && \mathcal{H}_p(\mathbb{C})
\end{tikzcd}
\]
Then we can write the operator $\Theta_p$ as
$$ \Theta_p:= \mathcal{E}_p \circ (\mathcal{V}_p)^{-1}.$$
By Theorem \ref{uni2} we know that the operator $ \mathcal{V}_p$ is unitary, hence its adjoint coincide with its inverse. So we have
$$ \Theta_p= \mathcal{E}_p \circ (\mathcal{V}_p)^{-1}=\mathcal{E}_p \circ \mathcal{E}_p=\left(\mathcal{I}+M_z \partial\right)^{2p}.$$
\end{proof}
\begin{theorem}
Let $p \in \mathbb{N}$. Then the inverse of the operator $\Theta_p$ is given by
$$ \Lambda_p=(\Theta_p)^{-1}=(R_0I)^{2p}.$$
Therefore, the operator $\Lambda_p$ maps the space $ \mathcal{H}_p(\mathbb{C})$ to $ \mathcal{F}_p(\mathbb{C})$.
\end{theorem}
\begin{proof}
By Theorem \ref{conn} we have that
$$ \Lambda_p=\left(\mathcal{E}_p \circ (\mathcal{V}_p)^{-1}\right)^{-1}= \mathcal{V}_p \circ (\mathcal{E}_p)^{-1}.$$
By Proposition \ref{adj3} we get
$$ \Lambda_p= \mathcal{V}_p \circ \mathcal{V}_p =(R_0I)^{2p}.$$
\end{proof}

By means of the operators $\Theta_p$ and $\Lambda_p$ we can realte the repducing kernel of the spaces $ \mathcal{F}_p(\mathbb{C})$ and $ \mathcal{H}_p(\mathbb{C})$

\begin{theorem}
\label{conn4}
Let $p \in \mathbb{N}$. Then for every $z, w \in \mathbb{C}$ we have
\begin{equation}
\label{conn0}
\Theta_p^{\bar{w}}\Theta_p^{z}(\mathcal{K}_p(z,w))=K_p(z,w),
\end{equation}
and
\begin{equation}
\label{conn1}
\Lambda_p^{\bar{w}}\Lambda_p^{z}(K_p(z,w))=\mathcal{K}_p(z,w),
\end{equation}
where we denote by $\Theta_p^{\bar{w}}$, $\Theta_p^{z}$, $\Lambda_p^{\bar{w}}$, $\Lambda_p^{z}$ the operators $\Theta_p$ and $\Lambda_p$ with respect the variables $z$ and $\bar{w}$, respectively.
\end{theorem}
\begin{proof}
We start by observing that $\Theta_p^{z}(z^n)=(n+1)^{2p} z^n$ and $\Theta_p^{\bar{w}}(\bar{w}^n)=(n+1)^{2p} \bar{w}^n$. Now by the definition of the kernel $\mathcal{K}_p(z,w)$ we have
\begin{eqnarray*}
\Theta_p^{\bar{w}}\Theta_p^{z}(\mathcal{K}_p(z,w))&=& \Theta_p^{\bar{w}}\Theta_p^{z}\left( \sum_{n=0}^{\infty} \frac{z^n \bar{w}^n}{(n+1)^{2p} n!}\right)\\
\nonumber
&=& \sum_{n=0}^{\infty} \frac{(n+1)^{2p} z^n \bar{w}^n}{ n!}\\
\nonumber
&=& K_p(z,w).
\end{eqnarray*}
The equality \eqref{conn1} follows from the fact that $\Lambda_p^{z}(z^n)=\frac{z^n}{(n+1)^{2p}} $, $\Lambda_p^{\bar{w}}(\bar{w}^n)=\frac{\bar{w}^n}{(n+1)^{2p}} $ and similar computations performed to get formula \eqref{conn0}.
\end{proof}

The above result implies the following connection between the generalized hypergeometric series and the Touchard polynomials.

\begin{corollary}
Let $p \in \mathbb{N}$. Then for every $z, w \in \mathbb{C}$ we have
$$ \Theta_p^{\bar{w}}\Theta_p^{z}({ }_{2p}F_{2p}(1,...,1; 2,...,2; z \bar{w}))=\frac{e^{z \bar{w}}}{z \bar{w}} T_{2p+1}(z\bar{w}),$$
and
$$ \Lambda_p^{\bar{w}}\Lambda_p^{z}\left(\frac{e^{z \bar{w}}}{z \bar{w}} T_{2p+1}(z\bar{w})\right)={ }_{2p}F_{2p}(1,...,1; 2,...,2; z \bar{w}),$$
\end{corollary}
\begin{proof}
The result follows by Theorem \ref{conn4}, Proposition \ref{rel} and Proposition \ref{Poch}.
\end{proof}

The operators $ \Theta_p$ and $\Lambda_p$ are also useful to relate the two type of Bargmann transforms introduced in this paper.

\begin{proposition}
Let $p \in \mathbb{N}$. Then we have
$$ \Theta_p( \mathcal{SB}_p )= \mathcal{B}_p,$$
and
$$ \Lambda_p(\mathcal{B}_p )=\mathcal{SB}_p.$$
\end{proposition}
\begin{proof}
By Proposition \ref{opBarg} we have
$$ \Theta_p (\mathcal{SB}_p)= \Theta_p \mathcal{V}_p (\mathcal{B}).$$
Since $ \Theta_{p}= \mathcal{E}_{2p}$ and by Theorem \ref{uni2} we have $ \mathcal{V}_{p}=(\mathcal{E}_p)^{-1}$, by using Proposition \ref{barp} we get
$$ \Theta_p (\mathcal{SB}_p)= \mathcal{E}_p(\mathcal{B} )= \mathcal{B}_p .$$
Now by Proposition \ref{barp} we have
$$\Lambda_p(\mathcal{B}_p )=\Lambda_p \mathcal{E}_p(\mathcal{B} ).$$
We observe that $\Lambda_p =\mathcal{V}_{2p}$ and by Proposition \ref{adj3} we have $\mathcal{E}_p=(\mathcal{V}_p)^{-1}$. Therefore by  Proposition \ref{opBarg} we get
$$\Lambda_p(\mathcal{B}_p )=\mathcal{V}_p(\mathcal{B} )=\mathcal{SB}_p .$$
\end{proof}

\section*{Acknowledgments} Daniel Alpay thanks the Foster G. and Mary McGaw Professorship in Mathematical Sciences, which supported this research. Kamal Diki thanks the Grand Challenges Initiative (GCI) at Chapman University, which supported this research.

\end{document}